\documentclass[10pt]{ijnam}
\input psfig.sty

\usepackage{stmaryrd}
\usepackage{amsmath}
\usepackage{amssymb}
\usepackage{tabularx,cases}
\usepackage{amsthm}
\usepackage{graphicx}
\usepackage{epstopdf}
\usepackage[bookmarksnumbered=true]{hyperref}
\def\currentvolume{x}
\def\currentissue{x}
\def\currentyear{xxxx}
\def\month{March}
\pagespan{1}{18}
\copyrightinfo{xxxx}{} % copyright year

\newcommand{\eps}{\varepsilon}
\newcommand{\bR}{\mathbb{R}}
\newcommand{\bx}{\mathbf{x}}

\newcommand{\fe}{\mathrm{e}}

\newcommand{\ud}{\mathrm{d}}

\newcommand{\be}{\begin{equation}}
\newcommand{\ee}{\end{equation}}
\newcommand{\ba}{\begin{array}}
\newcommand{\ea}{\end{array}}
\newcommand{\bea}{\begin{eqnarray}}
\newcommand{\eea}{\end{eqnarray}}
\newcommand{\beas}{\begin{eqnarray*}}
\newcommand{\eeas}{\end{eqnarray*}}
\newtheorem{remark}{Remark}[section]

\newtheorem{theorem}{Theorem}[section]

\numberwithin{equation}{section}
\begin{document}
%%%%%%%%%%%%%%%%

\title[High order Gautschi-type EWIs spectral method for the KGE]
{Symmetric high order Gautschi-type exponential wave integrators pseudospectral method for the nonlinear Klein-Gordon equation in the nonrelativistic limit regime}

% two authors have the same address
\author[Y. Wang]{Yan Wang}
\address{
  Beijing Computational Science Research Center, Beijing  100193, China
}
\email{matwyan@csrc.ac.cn}
\author[X. Zhao]{Xiaofei Zhao}
\address{
  IRMAR, Universit\'{e} de Rennes 1, France
}
\email{zhxfnus@gmail.com}
% Received by the editors ?
\date{}

% it is suggested to put it in Acknowledgments


\abstract{A group of high order Gautschi-type exponential wave integrators (EWIs) Fourier pseudospectral method are proposed and analyzed for solving the nonlinear Klein-Gordon equation (KGE) in the nonrelativistic limit regime, where a parameter $0<\eps\ll1$ which is inversely proportional to the speed of light, makes the solution propagate waves with wavelength $O(\eps^2)$ in time and $O(1)$ in space. With the Fourier pseudospectral method to discretize the KGE in space, we propose a group of EWIs with designed Gautschi's type quadratures for the temporal integrations, which can offer any intended even order of accuracy provided that the solution is smooth enough, while all the current existing EWIs offer at most second order accuracy. The scheme is explicit, time symmetric and rigorous error estimates show the meshing strategy of the proposed method is time step $\tau=O(\eps^2)$ and mesh size $h=O(1)$ as $0<\eps\ll1$, which is `optimal' among all classical numerical methods towards solving the KGE directly in the limit regime, and which also distinguish our methods from
other high order approaches such as Runge-Kutta methods which require $\tau=O(\eps^3)$. Numerical experiments with comparisons are done to confirm the error bound and show the superiority of the proposed methods over existing classical numerical methods.}


\maketitle

\section{Introduction}\label{sec: intro}
The Klein-Gordon equation (KGE) is known as the relativistic version of the Schr\"{o}dinger equation for describing the dynamics of spinless particles \cite{Sakurai}. Under proper nondimensionalization,
the dimensionless nonlinear KGE in $d$ dimensions $(d=1,2,3)$ reads
\cite{Dong,Machihara1,Machihara2,Masmoudi,Faou,Ginibre1,Ginibre2,Strauss}:
\begin{equation}\label{KG}
\left\{
  \begin{split}
    & \eps^2\partial_{tt}u-\Delta u+\frac{1}{\eps^2}u
    +f\left(u\right)=0,\quad \mathbf{x}\in\bR^d,\quad t>0,\\
    & u(\mathbf{x},0)= \phi_1(\mathbf{x}),\quad\partial_tu(\mathbf{x},0)
    =\frac{1}{\eps^2}\phi_2(\mathbf{x}),\quad \bx\in\bR^d.
  \end{split}
\right.
\end{equation}
Here $t$ is time, $\mathbf{x}$ is the spatial coordinate,
$u:=u(\mathbf{x},t)$ is a real-valued
scalar field, $0<\eps\leq1$ is a dimensionless parameter which is
inversely proportional to the speed of light,
 $\phi_1$ and $\phi_2$ are
two given real-valued initial data which are independent of $\eps$,
and $f(u): \bR\rightarrow\bR$ is a given nonlinearity independent of $\eps$.
It is clear that the KGE (\ref{KG}) is time symmetric and conserves the \emph{energy} \cite{Dong,Ginibre1,Ginibre2,Masmoudi}
\begin{align}\label{energy}
E(t)&:=\int_{\bR^d}\left[\eps^2|\partial_t u(\bx,t)|^2+|\nabla u(\bx,t)|^2+\frac{1}{\eps^2}|u(\bx,t)|^2+F(u(\bx,t))\right]d \bx\\
&\equiv\int_{\bR^d}\left[\frac{1}{\eps^2}|\phi_2(\bx)|^2+|\nabla \phi_1(\bx)|^2+\frac{1}{\eps^2}|\phi_1(\bx)|^2+F(\phi_1(\bx))\right]d \bx=E(0),\ t\ge0,\nonumber
\end{align}
with $F(u)=2\int_0^{u}f(\rho)d \rho$.

\begin{figure}[t!]
\centering
\includegraphics[height=6cm,width=12.5cm]{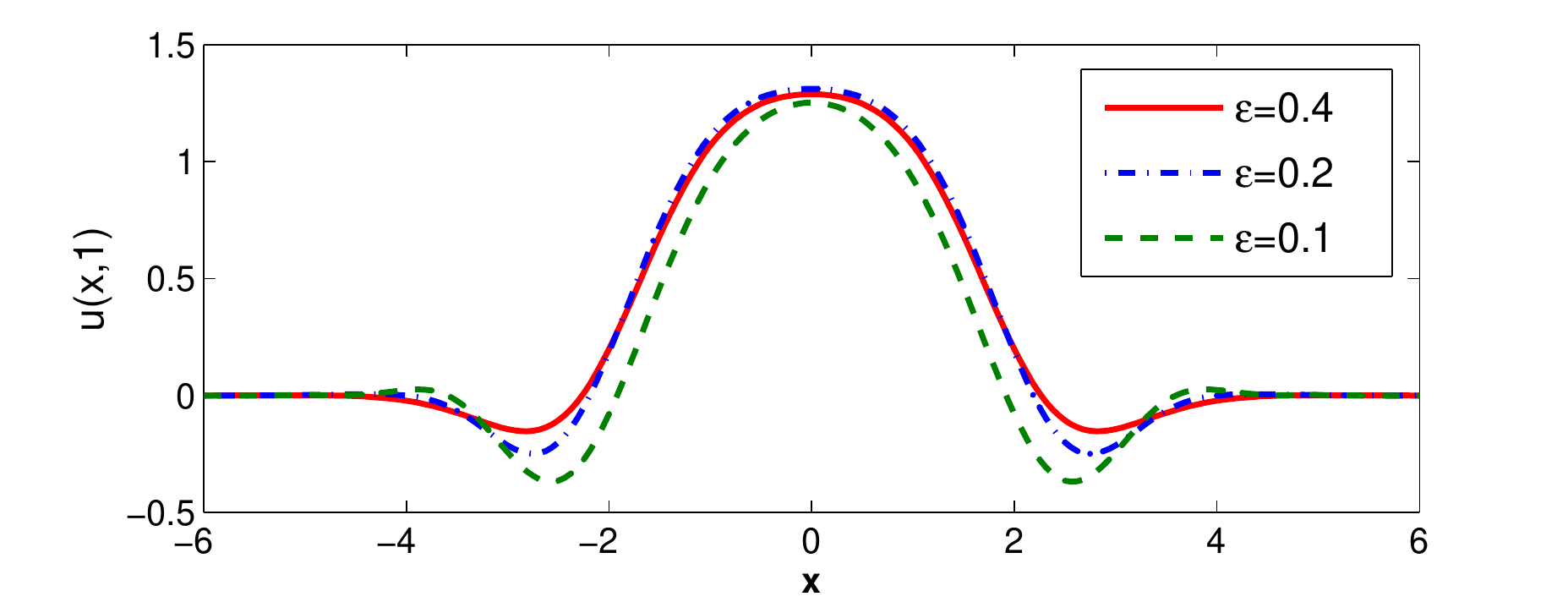}
\includegraphics[height=6cm,width=12.5cm]{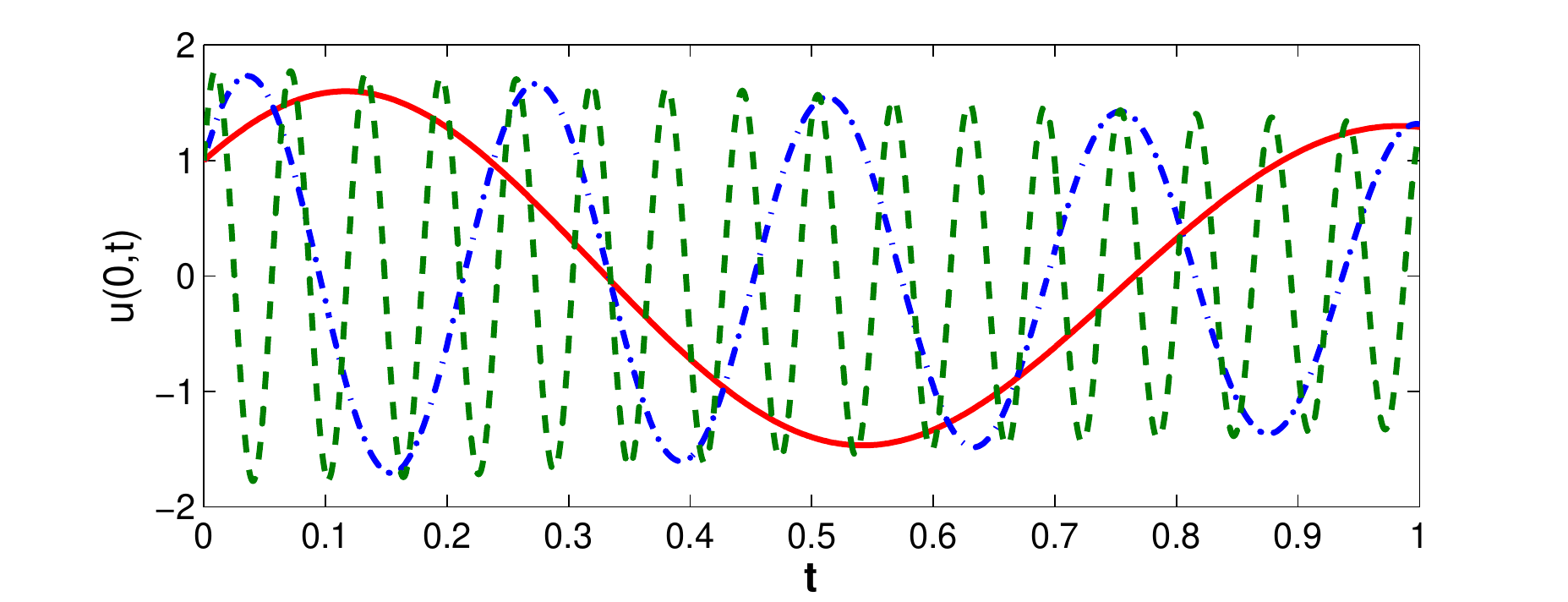}
\caption{The solution of (\ref{KG}) with $d=1$, $f(u)=u^3$, $\phi_1(x)=e^{-x^2/2}$ and
$\phi_2(x)=\frac{3}{2}\phi_1(x)$ for different $\eps$.}\label{fig:00}
\end{figure}

For fixed $0<\eps\leq1$, i.e. the relativistic regime, the KGE (\ref{KG}) has been well-studied both theoretically and numerically. We refer the readers to \cite{Dong} for a detailed review on the well-posedness and existing numerical methods for the KGE in this regime.
As $\eps\to0$, which corresponds to the speed of light goes to infinite and is known as the nonrelativistic limit in physics, recent analytical results \cite{Machihara1,Machihara2,Masmoudi} show that the problem (\ref{KG})
propagates waves with amplitude at $O(1)$, and wavelength
at $O(\eps^2)$ and $O(1)$ in time and space, respectively. The small wavelength makes the solution of the KGE highly oscillatory in time as $0<\eps\ll1$. Figure \ref{fig:00} shows an example of the profile of the solution under different $\eps$. The high oscillations cause severe numerical burdens in practical computations of the KGE in the nonrelativistic limit regime. For example, in order to capture the solution correctly in the highly oscillatory regime, frequently used finite difference time domain (FDTD) methods, such as the energy conservative type, semi-implicit type and fully explicit type \cite{Duncan,Strauss}, need the meshing strategy requirement (or $\eps$-scalability) $h=O(1)$ but $\tau=O(\eps^3)$ \cite{Dong}, where $h$ and$\tau$ denote the spatial mesh size and the time step, respectively. To release the temporal meshing strategy, based on the classical exponential wave integrators (EWIs) established in \cite{Hochbruck,Lubich1,Lubich2,Sanz,Gaustchi} for solving the oscillatory ODEs arising mainly from molecular dynamics, an EWI with the Gautschi-type quadrature \cite{Gaustchi} spectral method was proposed for solving the nonlinear KGE in the nonrelativistic limit regime and was shown to improve the $\eps$-scalability to $\tau=O(\eps^2)$ in \cite{Dong}. This method also finds successful applications in solving the Klein-Gordon-Zakharov (KGZ) system in a similar oscillatory situation \cite{BDZ}. Later on, an EWI with the Deuflhard-type quadrature \cite{Deuflhard} spectral method, which is equivalent to the time-splitting spectral method, was considered in \cite{ZhaoCicp} for the KGE in the nonrelativistic limit regime. It can offer a smaller temporal error bounded but the same $\eps$-scalability. Recent studies turn to utilize multiscale analysis to first derive some sophisticate reformulations or decompositions of the KGE, then based on which one can propose some suitable numerical methods \cite{Faou,Chartier,ZhaoMTI,Zhao} for asymptotic preserving or uniformly accurate property. These multiscale numerical methods are extremely powerful in computations of KGE in the oscillatory regime, however they either require some delicate pre-knowledge of the oscillation structures of the problem \cite{ZhaoMTI,Faou,Zhao} or require introducing an extra degree-of-freedom \cite{Chartier}. In view of that the solution to (\ref{KG}) has oscillation wavelength at $O(\eps^2)$ in time, the EWIs could be viewed as the optimal one among all the traditional methods towards
integrating the KGE (\ref{KG}) directly in the nonrelativistic limit regime.

However, all the existing EWIs for either solving the oscillatory ODEs from molecular dynamics or solving the KGE offer at most second order accuracy in temporal discretization. Of course, one can apply the Runge-Kutta methods, like the one proposed in \cite{Dehghan} for the approximations in time to get higher order temporal convergence rates, but that will surely lead to lost of the time symmetry property or worse $\eps$-scalability in the nonrelativistic limit regime. The time symmetry is known as a key property to provide good long time behaviors of the numerical schemes \cite{Lubich1,Lubich2}. This work is devoted to propose a group of symmetric Gautschi-type EWIs with high order convergence rate in temporal approximation and with Fourier pseudospectral discretization in space for solving the KGE (\ref{KG}) in the nonrelativistic limit regime. We are going to apply the Fourier spectral method for the spatial discretization of the KGE at first, then propose a group of symmetric Gautschi-type EWIs with even order of accuracy for integrating the oscillatory ODEs resulting from the spatial semi-discretization in the Fourier frequency space. The scheme is fully explicit, easy to implement, and in principle, one can construct the scheme to get any even order of temporal accuracy provided that the solution of the KGE (\ref{KG}) is smooth enough. Rigorous error estimates of the proposed methods are established in the highly oscillatory regime with particular attentions paid to the dependence of $\eps$ in the error bound, where the results show that the $\eps$-scalability of the high order Gautschi-type EWIs spectral method is still $h=O(1)$ and $\tau=O(\eps^2)$ as $0<\eps\ll1$. It is believed that the proof techniques here could also give some clues to the error estimates of a group of trigonometric integrators considered in \cite{DongNew}. It is also believed that the higher order EWIs could offer an essential way to extend the order of the uniform accuracy of the recent developed multiscale time integrators \cite{Zhao,ZhaoMTI,CaiDirac,ZhaoKGS}. Extensive numerical experiments will justify the theoretical error bound and comparisons with the existing classical numerical methods will show the superiority of the high order methods in accuracy and energy preserving property. The proposed high order EWIs spectral method could also find applications to solve other KG-type equations or coupled system, such as the KGZ system or the Klein-Gordon-Schr\"{o}dinger system \cite{BaoYang,ZhaoKGS}.

The rest of the paper is organized as follows. In Section \ref{sec: method}, we derive the high order EWIs spectral method. The convergence theorem with rigorous proof is given in Section \ref{sec: thm}. Numerical results are reported in Section \ref{sec: result}. Finally, some concluding remarks are drawn in Section \ref{sec: con}. Throughout this paper,
we adopt the notation $A\lesssim B$ to represent that there exists a generic constant $C>0$,
which is independent of time step $\tau$ (or $n$), mesh size $h$ and $\eps$, such that $|A|\leq CB$.

\section{Numerical methods}\label{sec: method}
In this section, we shall first derive a detailed 4th order Gautschi-type EWI Fourier pseudospectral method for solving the KGE (\ref{KG}), and then present the general even order EWIs with spectral discretization. For the simplicity of notations, we present the numerical method
in one space dimension (1D), i.e. $d=1$ in (\ref{KG}). Generalizations to higher
dimensions are straightforward and results remain valid without modifications.
Due to fast decay of the solution of the KGE (\ref{KG}) at far field \cite{Machihara1,Machihara2,Masmoudi}, similar to those in the literature for numerical computations \cite{Dong,Duncan,Faou,Jimenez,Strauss},
the whole space problem (\ref{KG}) in 1D is truncated
onto a finite interval $\Omega=(a,b)$ with periodic boundary conditions
($a$ and $b$ are usually chosen sufficient large such that
the truncation error is negligible):
\begin{numcases}
\ \eps^2\partial_{tt} u(x,t)-\partial_{xx} u(x,t)+\frac{1}{\eps^2}u(x,t)
+f\left(u(x,t)\right)=0,\quad x\in\Omega,\ t>0,\nonumber\\
 u(a,t)=u(b,t),\quad \partial_x u(a,t)=\partial_x u(b,t), \qquad t\geq0,\label{KG trun}\\
 u(x,0)= \phi_1(x),\quad\partial_tu(x,0)=\frac{1}{\eps^2}\phi_2(x),\qquad x\in\overline{\Omega}=[a,b].\nonumber
\end{numcases}

\subsection{4th order Gautschi-type EWI}
Choose the mesh size $h:=\Delta x=(b-a)/M$ with $M$ a positive even integer
 and denote grid points as
$x_j:=a+jh$ for  $j=0,1,\ldots, M$. Define
\begin{align*}
&X_M:=\mbox{span}\left\{\psi_l(x)=\fe^{i\mu_l(x-a)}\ |\ \mu_l=\frac{2\pi l}{b-a},\ l=-\frac{M}{2}, \ldots, \frac{M}{2}-1\right\},
 \\
&Y_M:=\left\{\mathbf{v}=(v_0,v_1,\ldots,v_M)\in\bR^{M+1} \; |\; v_0=v_M\right\}\quad \hbox{with} \
 \|\mathbf{v}\|_{l^2}^2=h\sum_{j=0}^{M-1}|v_j|^2.
\end{align*}
For a periodic function $v(x)$ on $\overline{\Omega}$ and a vector $\mathbf{v}\in Y_M$,
let $P_M: L^2(\Omega)\rightarrow X_M$ be the standard $L^2$-projection operator,
 and $I_M: C(\Omega)\rightarrow X_M$ or  $Y_M \rightarrow X_M$
be the trigonometric interpolation operator \cite{Shen,book-sp-GO,book-sp-HGG}, i.e.
\begin{equation}\label{project oper}
(P_Mv)(x)=\sum_{l=-M/2}^{M/2-1}\widehat{v}_l\;\psi_l(x),\quad (I_M\mathbf{v})(x)=\sum_{l=-M/2}^{M/2-1}\widetilde{v}_l\;\psi_l(x),\quad a\leq x\leq b,
\end{equation}
where $\widehat{v}_l$ and $\widetilde{v}_l$ are the Fourier and discrete Fourier transform
coefficients of the periodic function $v(x)$ and vector $\mathbf{v}$, respectively, defined as
\begin{equation}\label{Fourier tans}
\widehat{v}_l=\frac{1}{b-a}\int_a^b v(x)\;\fe^{-i\mu_l(x-a)}dx,\qquad \widetilde{v}_l=\frac{1}{M}\sum_{j=0}^{M-1}v_j\;\fe^{-i\mu_l(x_j-a)}.
\end{equation}
Let $\tau=\Delta t>0$ be the step size, and denote time steps by $t_n=n\tau$ for $n=0,1,\ldots$.
Then a Fourier spectral method for discretizing (\ref{KG trun}) with $t=t_n+s\,(n=0,1,\ldots)$ reads:
Find $u_{M}(x,t_n+s)\in X_M$, i.e.
\begin{equation} \label{uM def}
u_{M}(x,t_n+s)=\sum_{l=-M/2}^{M/2-1}\widehat{(u_M)}_l(t_n+s)\psi_l(x),\quad x\in\Omega,\ s\in\bR,
\end{equation}
such that
\begin{equation}\label{uM KG}
\eps^2\partial_{ss} u_M(x,t_n+s)-\partial_{xx} u_M(x,t_n+s)+\frac{1}{\eps^2}u_M(x,t_n+s)
+P_Mf\left(u_M(x,t_n+s)\right)=0.
\end{equation}
Plugging (\ref{uM def}) into (\ref{uM KG}), and noticing the orthogonality of $\psi_l(x)$ for $l=-\frac{M}{2},\ldots,\frac{M}{2}-1$, we get for $n=0,1,\ldots,$
\begin{equation}\label{KG Fourier}
\eps^2\widehat{(u_M)}_l''(t_n+s)+\left(\mu_l^2+\frac{1}{\eps^2}\right)\widehat{(u_M)}_l(t_n+s)
+\widehat{(f_M^n)}_l(s)=0,\quad s\in\bR,
\end{equation}
where $f_M^n(x,s):=P_Mf\left(u_M(x,t_n+s)\right).$
By using the variation-of-constant formula to (\ref{KG Fourier}), we get for $n=0,1,\ldots,$
\begin{align}\label{vcf s}
\widehat{(u_M)}_l(t_n+s)=&\cos(\omega_l s)\widehat{(u_M)}_l(t_n)+\frac{\sin(\omega_l s)}{\omega_l}\widehat{(u_M)}_l'(t_n)\\
&-\int_0^s\frac{\sin(\omega_l(s-w))}{\eps^2\omega_l}\widehat{(f_M^n)}_l(w)\ud w,\quad  l=-\frac{M}{2},\ldots,\frac{M}{2}-1,\nonumber
\end{align}
with $\omega_l:=\frac{\sqrt{\eps^2\mu_l^2+1}}{\eps^2}$.
Differentiating (\ref{vcf s}) with respect to $s$ on both sides, we get
\begin{align}\label{vcf ds}
\widehat{(u_M)}_l'(t_n+s)=&-\omega_l\sin(\omega_l s)\widehat{(u_M)}_l(t_n)+\cos(\omega_l s)\widehat{(u_M)}_l'(t_n)\\
&-\int_0^s\frac{\cos(\omega_l(s-w))}{\eps^2}\widehat{(f_M^n)}_l(w)\ud w,\quad l=-\frac{M}{2},\ldots,\frac{M}{2}-1.\nonumber
\end{align}
For $n\geq1$, changing $s$ to $-s$ in (\ref{vcf s}) and (\ref{vcf ds}), we get
\begin{subequations}\label{vcf -s}
\begin{align}
\widehat{(u_M)}_l(t_n-s)=&\cos(\omega_l s)\widehat{(u_M)}_l(t_n)-\frac{\sin(\omega_l s)}{\omega_l}\widehat{(u_M)}_l'(t_n)\label{vcf -sa}\\
&-\int_0^s\frac{\sin(\omega_l(s-w))}{\eps^2\omega_l}\widehat{(f_M^n)}_l(-w)\ud w,\nonumber\\
\widehat{(u_M)}_l'(t_n-s)=&\omega_l\sin(\omega_l s)\widehat{(u_M)}_l(t_n)+\cos(\omega_l s)\widehat{(u_M)}_l'(t_n)\label{vcf -sb}\\
&+\int_0^s\frac{\cos(\omega_l(s-w))}{\eps^2}\widehat{(f_M^n)}_l(-w)\ud w,\quad l=-\frac{M}{2},\ldots,\frac{M}{2}-1.\nonumber
\end{align}
\end{subequations}
Adding (\ref{vcf -sa}) to (\ref{vcf s}) and subtracting (\ref{vcf -sb}) from (\ref{vcf ds}) for $n\geq1$, and then let $s=\tau$, we get
\begin{subequations}\label{vcf}
\begin{align}
\widehat{(u_M)}_l(t_{n+1})=&-\widehat{(u_M)}_l(t_{n-1})+2\cos(\omega_l\tau)\widehat{(u_M)}_l(t_n)\\
&-\int_0^\tau\frac{\sin(\omega_l(\tau-w))}{\eps^2\omega_l}
\left[\widehat{(f_M^n)}_l(w)+\widehat{(f_M^n)}_l(-w)\right]\ud w,\nonumber\\
\widehat{(u_M)}_l'(t_{n+1})=&\widehat{(u_M)}_l'(t_{n-1})-2\omega_l\sin(\omega_l \tau)\widehat{(u_M)}_l(t_n)\\
&-\int_0^\tau\frac{\cos(\omega_l(\tau-w))}{\eps^2}\left[\widehat{(f_M^n)}_l(w)+\widehat{(f_M^n)}_l(-w)\right]\ud w.\nonumber
\end{align}
\end{subequations}
Similar to the quadrature proposed by W. Gaustchi in \cite{Gaustchi} and used in \cite{Dong}, but in order to construct a fourth order accuracy method,  we approximate the unknown integrals in (\ref{vcf}) by using the Taylor's expansion of the nonlinearity up to the second order terms as
\begin{subequations}\label{quadrature}
\begin{align}
&\quad\int_0^\tau\frac{\sin(\omega_l(\tau-w))}{\eps^2\omega_l}\left[\widehat{(f_M^n)}_l(w)+\widehat{(f_M^n)}_l(-w)\right]\ud w\nonumber\\
&\approx
\int_0^\tau\frac{\sin(\omega_l(\tau-w))}{\eps^2\omega_l}\left[2\widehat{(f_M^n)}_l(0)+w^2\frac{\ud^2}{\ud s^2}\widehat{(f_M^n)}_l(0)\right]\ud w,\\
&\quad\int_0^\tau\frac{\cos(\omega_l(\tau-w))}{\eps^2}\left[\widehat{(f_M^n)}_l(w)+\widehat{(f_M^n)}_l(-w)\right]\ud w\nonumber\\
&\approx
\int_0^\tau\frac{\cos(\omega_l(\tau-w))}{\eps^2}\left[2\widehat{(f_M^n)}_l(0)+w^2\frac{\ud^2}{\ud s^2}\widehat{(f_M^n)}_l(0)\right]\ud w,\quad n\geq1,
\end{align}
\end{subequations}
and then carry out the rest trigonometric integrations exactly, where
\begin{equation}\label{coefficient}
\left\{
\begin{split}
&\int_0^\tau\frac{\sin(\omega_l(\tau-w))}{\eps^2\omega_l}\ud w=\frac{1}{\eps^2\omega_l^2}[1-\cos(\omega_l\tau )],\\
&\int_0^\tau\frac{\sin(\omega_l(\tau-w))}{\eps^2\omega_l}w^2\ud w=\frac{1}{\eps^2\omega_l^4}
\left[\omega_l^2\tau^2+2\cos(\omega_l\tau)-2\right],\\
&\int_0^\tau\frac{\cos(\omega_l(\tau-w))}{\eps^2}\ud w=\frac{1}{\eps^2\omega_l}\sin(\omega_l\tau ),\\
&\int_0^\tau\frac{\cos(\omega_l(\tau-w))}{\eps^2}w^2\ud w=\frac{1}{\eps^2\omega_l^3}
\left[2\omega_l\tau-2\sin(\omega_l\tau)\right].
\end{split}\right.
\end{equation}
For the second order derivatives involved in the above approximations (\ref{quadrature}), i.e.
$$\frac{\ud^2}{\ud s^2}\widehat{(f_M^n)}_l(0)=\widehat{\left(\partial_{ss}f(u_M(x,t_n+s))\right)}_l\big|_{s=0},\quad n\geq1,$$ it can be found out from the equation (\ref{uM KG}), i.e.
\begin{align*}
&\quad\partial_{ss}f(u_M(x,t_n+s))\big|_{s=0}\\
&=f'(u_M(x,t_n))\partial_{ss}u_M(x,t_n)+f''(u_M(x,t_n))(\partial_{s}u_M(x,t_n))^2\\
&=\frac{1}{\eps^2}f'(u_M(x,t_n))\left[\partial_{xx} u_M(x,t_n)-\frac{1}{\eps^2}u_M(x,t_n)-P_Mf(u_M(x,t_n))\right]\\
&\quad\ +f''(u_M(x,t_n))(\partial_{s}u_M(x,t_n))^2,\quad n\geq1.
\end{align*}
Since the numerical integrator based on (\ref{vcf}) proceeds in a three-level format, we need to find approximations of $u_M(x,t_1)$ and $\partial_su_M(x,t_1)$ to start the scheme. Taking $n=0$ and $s=\tau$ in (\ref{vcf s}) and (\ref{vcf ds}), we get
\begin{subequations}\label{vcf 1}
\begin{align}
 \widehat{(u_M)}_l(t_1)=&\cos(\omega_l\tau)\widehat{(u_M)}_l(0)+\frac{\sin(\omega_l\tau)}{\omega_l}\widehat{(u_M)}_l'(0)\\
&-\int_0^\tau\frac{\sin(\omega_l(\tau-w))}{\eps^2\omega_l}\widehat{(f_M^0)}_l(w)\ud w,\nonumber\\
\widehat{(u_M)}_l'(t_1)=&-\omega_l\sin(\omega_l\tau)\widehat{(u_M)}_l(0)+\cos(\omega_l s)\widehat{(u_M)}_l'(0)\\
&-\int_0^\tau\frac{\cos(\omega_l(\tau-w))}{\eps^2}\widehat{(f_M^0)}_l(w)\ud w,\qquad l=-\frac{M}{2},\ldots,\frac{M}{2}-1.\nonumber
\end{align}
\end{subequations}
Applying the quadrature similar as (\ref{quadrature}) to the unknown integrals in (\ref{vcf 1}), i.e.
\begin{subequations}
\begin{align*}
&\quad\int_0^\tau\frac{\sin(\omega_l(\tau-w))}{\eps^2\omega_l}\widehat{(f_M^0)}_l(w)\ud w\\
&\approx\int_0^\tau\frac{\sin(\omega_l(\tau-w))}{\eps^2\omega_l}\left[\widehat{(f_M^0)}_l(0)+w\frac{\ud}{\ud s}\widehat{(f_M^0)}_l(0)+\frac{w^2}{2}\frac{\ud^2}{\ud s^2}\widehat{(f_M^0)}_l(0)\right]\ud w,\\
&\quad\int_0^\tau\frac{\cos(\omega_l(\tau-w))}{\eps^2}\widehat{(f_M^0)}_l(w)\ud w\\
&\approx\int_0^\tau\frac{\cos(\omega_l(\tau-w))}{\eps^2}\left[\widehat{(f_M^0)}_l(0)+w\frac{\ud}{\ud s}\widehat{(f_M^0)}_l(0)+\frac{w^2}{2}\frac{\ud^2}{\ud s^2}\widehat{(f_M^0)}_l(0)\right]\ud w,
\end{align*}
\end{subequations}
where we have $\frac{\ud}{\ud s}\widehat{(f_M^0)}_l(0)=\widehat{\left(\partial_{s}f(u_M(x,s))\right)}_l\big|_{s=0}$, and in addition to (\ref{coefficient}),
\begin{equation}\label{coefficient0}
\left\{
\begin{split}
&\int_0^\tau\frac{\sin(\omega_l(\tau-w))}{\eps^2\omega_l}w\ud w=\frac{1}{\eps^2\omega_l^3}\left[\omega_l\tau-\sin(\omega_l\tau)\right],\\
&\int_0^\tau\frac{\cos(\omega_l(\tau-w))}{\eps^2}w\ud w=\frac{1}{\eps^2\omega_l^2}\left[1-\cos(\omega_l\tau)\right].
\end{split}\right.
\end{equation}
The above temporal approximations for both (\ref{vcf}) and (\ref{vcf 1}) offer naturally fourth order truncation error bounds and clearly become exact when the nonlinearity $f(\cdot)$ reduces to a constant function.

In details, a 4th order Gautschi-type EWI Fourier spectral method (4th-GIFS) reads as follows.
Denote $u_M^n(x)$ and $\dot{u}_M^n(x)\,(n=0,1,\ldots)$ be the approximations to $u(x,t_n)$ and $\partial_tu(x,t_n)$, respectively. Choose $u_M^0(x)=\phi_1(x)$ and $\dot{u}_M^0(x)=\frac{1}{\eps^2}\phi_2(x)$, then for $n\geq0,$
\begin{equation}\label{4th GIFS s}
u_{M}^{n+1}(x)=\sum_{l=-M/2}^{M/2-1}\widehat{(u_M^{n+1})}_l\psi_l(x),\quad
\dot{u}_{M}^{n+1}(x)=\sum_{l=-M/2}^{M/2-1}\widehat{(\dot{u}_M^{n+1})}_l\psi_l(x),\quad x\in\Omega,
\end{equation}
where
\begin{subequations}\label{4th GIFS m}
\begin{align}
\widehat{(u_M^{n+1})}_l=&-\widehat{(u_M^{n-1})}_l+2\cos(\omega_l\tau)\widehat{(u_M^{n})}_l-\frac{2-2\cos(\omega_l\tau)}{\eps^2\omega_l^2}
\widehat{f^n_l}\nonumber\\
&-\frac{\omega_l^2\tau^2+2\cos(\omega_l\tau)-2}{\eps^2\omega_l^4}\widehat{\ddot{f}^n_l},\qquad n\geq1,\\
\widehat{(\dot{u}_M^{n+1})}_l=&\widehat{(\dot{u}_M^{n-1})}_l-2\omega_l\sin(\omega_l \tau)\widehat{(\dot{u}_M^{n})}_l-
\frac{2\sin(\omega_l\tau)}{\eps^2\omega_l}\widehat{f^n_l}\nonumber\\
&-\frac{2\omega_l\tau -2\sin(\omega_l\tau)}{\eps^2\omega_l^3}\widehat{\ddot{f}^n_l},\qquad n\geq1,
\end{align}
\end{subequations}
and
\begin{subequations}\label{4th GIFS e}
\begin{align}
\widehat{(u_M^{1})}_l=&\cos(\omega_l\tau)\widehat{(\phi_1)}_l+\frac{\sin(\omega_l \tau)}{\eps^2\omega_l}\widehat{(\phi_2)}_l-\frac{1-\cos(\omega_l\tau)}{\eps^2\omega_l^2}\widehat{f^0_l}\nonumber\\
&-\frac{\omega_l\tau-\sin(\omega_l\tau)}{\eps^2\omega_l^3}\widehat{\dot{f}^0_l}-
\frac{\omega_l^2\tau^2+2\cos(\omega_l\tau)- 2}{2\eps^2\omega_l^4}\widehat{\ddot{f}^0_l},\\
\widehat{(\dot{u}_M^{1})}_l=&-\omega_l\sin(\omega_l\tau)\widehat{(\phi_1)}_l+\frac{\cos(\omega_l \tau)}{\eps^2}\widehat{(\phi_2)}_l-\frac{\sin(\omega_l\tau)}{\eps^2\omega_l}\widehat{f^0_l}\nonumber\\
&-\frac{1-\cos(\omega_l\tau)}{\eps^2\omega_l^2}\widehat{\dot{f}^0_l}-\frac{\omega_l\tau-\sin(\omega_l\tau)}{\eps^2\omega_l^3}\widehat{\ddot{f}^0_l},
\end{align}
\end{subequations}
with
\begin{align*}
&f^n(x)=f(u_M^n(x)),\qquad \dot{f}^0(x)=\frac{1}{\eps^2}f'\left(\phi_1(x)\right)\cdot\phi_2(x),\qquad n\geq0,\\
&\ddot{f}^n(x)=f''(u_M^n(x))\cdot\left(\dot{u}_M^n(x)\right)^2+\frac{1}{\eps^2}f'(u_M^n(x))
\cdot\left[\partial_{xx}u_M^n(x)-\frac{1}{\eps^2}u_M^n(x)-f^n(x)\right].
\end{align*}
In practice, the integrals defined in (\ref{Fourier tans}) for computing the Fourier transform
coefficients in (\ref{4th GIFS s})-(\ref{4th GIFS e}) are not suitable, and they are usually replaced by the interpolations as defined in (\ref{Fourier tans}) \cite{Shen,Dong,Gottlieb}.
Thus, a 4th order Gautschi-type EWI Fourier pseudospectral method (4th-GIFP) reads as follows. Let $u_j^n$ and $\dot{u}_j^n\,(n=0,1,\ldots,j=0,\ldots,M)$ be the approximations to $u(x_j,t_n)$ and $\partial_tu(x_j,t_n)$, respectively. Choose $u_j^0=\phi_1(x_j)$ and $\dot{u}_j^0=\frac{1}{\eps^2}\phi_2(x_j)$, then for $n\geq0$,
\begin{equation}\label{4th GIFP s}
u_{j}^{n+1}=\sum_{l=-M/2}^{M/2-1}\widetilde{u^{n+1}_l}\psi_l(x_j),\quad
\dot{u}_{j}^{n+1}=\sum_{l=-M/2}^{M/2-1}\widetilde{\dot{u}^{n+1}_l}\psi_l(x_j),\quad j=0,\ldots,M,
\end{equation}
where
\begin{subequations}\label{4th GIFP m}
\begin{align}
\widetilde{u^{n+1}_l}=&-\widetilde{u^{n-1}_l}+2\cos(\omega_l\tau)\widetilde{u^{n}_l}-\frac{2-2\cos(\omega_l\tau)}{\eps^2\omega_l^2}
\widetilde{f^n_l}+\frac{2-2\cos(\omega_l\tau)-\omega_l^2\tau^2}{\eps^2\omega_l^4}\widetilde{\ddot{f}^n_l},\\
\widetilde{\dot{u}^{n+1}_l}=&\widetilde{\dot{u}^{n-1}_l}-2\omega_l\sin(\omega_l \tau)\widetilde{\dot{u}^{n}_l}-
\frac{2\sin(\omega_l\tau)}{\eps^2\omega_l}\widetilde{f^n_l}+\frac{2\sin(\omega_l\tau)-2\omega_l\tau}{\eps^2\omega_l^3}\widetilde{\ddot{f}^n_l},\quad n\geq1,
\end{align}
\end{subequations}
and
\begin{subequations}\label{4th GIFP e}
\begin{align}
\widetilde{u^{1}_l}=&\cos(\omega_l\tau)\widetilde{(\phi_1)}_l+\frac{\sin(\omega_l \tau)}{\eps^2\omega_l}\widetilde{(\phi_2)}_l-\frac{1-\cos(\omega_l\tau)}{\eps^2\omega_l^2}\widetilde{f^0_l}\nonumber\\
&-\frac{\omega_l\tau-\sin(\omega_l\tau)}{\eps^2\omega_l^3}\widetilde{\dot{f}^0_l}-\frac{\omega_l^2\tau^2+2\cos(\omega_l\tau)- 2}{2\eps^2\omega_l^4}\widetilde{\ddot{f}^0_l},\\
\widetilde{\dot{u}^{1}_l}=&-\omega_l\sin(\omega_l\tau)\widetilde{(\phi_1)}_l+\frac{\cos(\omega_l \tau)}{\eps^2}\widetilde{(\phi_2)}_l-\frac{\sin(\omega_l\tau)}{\eps^2\omega_l}\widetilde{f^0_l}\nonumber\\
&-\frac{1-\cos(\omega_l\tau)}{\eps^2\omega_l^2}\widetilde{\dot{f}^0_l}-\frac{\omega_l\tau-\sin(\omega_l\tau)}{\eps^2\omega_l^3}\widetilde{\ddot{f}^0_l},
\end{align}
\end{subequations}
with
\begin{align*}
&f^n_j=f(u_j^n),\qquad \dot{f}^0_j=\frac{1}{\eps^2}f'\left(\phi_1(x_j)\right)\cdot\phi_2(x_j),\qquad n\geq0,\\
&\ddot{f}^n_j=f''(u_j^n)\cdot\left(\dot{u}_j^n\right)^2+\frac{1}{\eps^2}f'(u_j^n)
\cdot\left[\partial_{xx}I_M(u^n)(x_j)-\frac{1}{\eps^2}u_j^n-f^n_j\right].
\end{align*}
Clearly, the proposed 4th-GIFS (\ref{4th GIFS s})-(\ref{4th GIFS e}) or 4th-GIFP (\ref{4th GIFP s})-(\ref{4th GIFP e}) is fully explicit and easy to implement. It is very efficient due to the fast Fourier transform (FFT), and its memory
cost is $O(M)$ and the computational cost per time step is $O(M \log M )$. The scheme is also clearly time symmetric, i.e. exchanging $n+1$ with $n-1$ and changing $\tau$ to $-\tau$ in (\ref{4th GIFS m}) or (\ref{4th GIFP m}), it remains the same.

\subsection{Higher order Gautschi-type EWIs spectral method}\label{subsec: EWIs}
As a natural generalization, one can approximate the nonlinearity in (\ref{vcf}) by using its Taylor's expansion up to some higher order like 6th, 8th...terms, provided that the nonlinearity and the solution to (\ref{KG}) are smooth enough. In this case, we can get an arbitrary $2N$th order accurate Gautschi-type EWI Fourier spectral/pesudospectral method (2Nth-GIFS/2Nth-GIFP) for some integer $N\geq2$ by using
\begin{align*}
\widehat{(f_M^n)}_l(w)+\widehat{(f_M^n)}_l(-w)=&2\widehat{(f_M^n)}_l(0)+2\sum_{m=1}^{N-1}\frac{w^{2m}}{(2m)!}\frac{\ud^{2m}}{\ud w^{2m}}\widehat{(f_M^n)}_l(0)\\
&+\mathcal{O}(w^{2N}),\quad 0\leq w\leq\tau,\ n\geq1,
\end{align*}
and then similar as before, carrying out the trigonometric integrations left in (\ref{vcf}) exactly. Also, for the starting values (\ref{vcf 1}), take the Taylor's expansion
$$\widehat{(f_M^0)}_l(w)=\sum_{m=0}^{2N-2}\frac{w^{m}}{m!}\frac{\ud^{m}}{\ud w^{m}}\widehat{(f_M^0)}_l(0)+\mathcal{O}(w^{2N-1}),\quad 0\leq w\leq\tau,$$
for the nonlinearity in (\ref{vcf 1}) in order to get a $2N$th order approximation and then evaluate the integrals. Consequently, we will need the higher order time derivatives of $u(x,t)$. This can be obtained from the original problem (\ref{KG}) with lower order derivatives in hands, i.e.
$$\partial_t^mu(x,t)=\frac{1}{\eps^2}\partial_t^{m-2}\left(\partial_{xx} u(x,t)-\frac{1}{\eps^2}u(x,t)-f(u(x,t))\right),\quad m\geq2.$$
Then the scheme of 2Nth-GIFS or 2Nth-GIFP can be written down similarly as (\ref{4th GIFS s})-(\ref{4th GIFS s}) or (\ref{4th GIFP s})-(\ref{4th GIFP s}). We omit the details here for brevity.

Again, the proposed high order GIFS/GIFP methods are fully explicit, time symmetric, efficient due to the FFT and become exact when $f(\cdot)$ is a constant.

To close this chapter, we make some remarks on the proposed method. The 4th order or higher order symmetric Gautschi-type EWIs Fourier spectral/pseudospectral method can be easily applied and extended to solve other KG-type equations or coupled system, such as the Klein-Gordon-Zakharov system in the high-plasma-frequency and subsonic limit regime \cite{BDZ} where similar oscillations occur. Similar numerical schemes with similar expected numerical performance can be derived. We also remark that if the periodic boundary condition for the KGE (\ref{KG trun}) is replaced by the
homogeneous Dirichlet or Neumann boundary condition which is also suitable here for domain truncations, the GIFS/GIFP method and its following error estimates
are still valid with the Fourier basis is replaced by sine or cosine basis. For some inhomogeneous general boundary conditions, one can turn to the compact finite difference discretization in order to get high order spatial accuracy and then construct similar high order Gautschi-type EWIs to the resulting ODEs from corresponding semidisretizations.

\section{Convergence result}\label{sec: thm}
In this section, we present the rigorous error estimate results of the proposed 4th order Gautschi-type EWI Fourier spectral/pseudospectral method for solving the KGE (\ref{KG trun}). Generalizations of the results to higher order EWIs spectral/pseudospectral method can be obtained similarly with stronger regularity assumptions on the solution.

\subsection{Main result}
To get the optimal error estimates for the 4th order scheme, we make assumptions on the solution of the KGE (\ref{KG trun}) motivated from \cite{Machihara1,Machihara2,Masmoudi} as:
\begin{equation}\tag{A}\label{assumption}
\begin{split}
&f(\cdot)\in C^2(\bR),\quad u\in C^1\left([0,T];H_p^{m_0+1}(\Omega)\right)\cap C^4\left([0,T];H^{1}(\Omega)\right),\\
&\|\partial_t^ku\|_{L^\infty([0,T];H^{m_0+1})}\lesssim\frac{1}{\eps^{2k}},\quad k=0,1;\quad\|\partial_t^ku\|_{L^\infty([0,T];H^{1})}\lesssim\frac{1}{\eps^{2k}},\quad k=3,4,
\end{split}
\end{equation}
where $m_0\geq2$, $0<T\leq T^*$ with $T^*$ the maximum existence time of the solution and
$$H_p^{m_0+1}(\Omega):=\left\{v\in H^{m_0+1}(\Omega):\,\partial_x^{m}v(a)=\partial_x^{m}v(b),\ m=0,\ldots,m_0\right\}.$$
Under assumption (A), denote
$$C_0=\max_{0<\eps\le 1}\left\{\left\|u\right\|_{L^\infty([0,T];H^1\cap L^{\infty})},
\ \eps^2\left\|\partial_tu\right\|_{L^\infty([0,T];H^1\cap L^{\infty})}\right\},$$
and with $u_M^n,\dot{u}_M^n$ obtained from the 4th-GIFS (\ref{4th GIFS s})-(\ref{4th GIFS e}), define the error functions as
\begin{equation}\label{error fun}
e^n(x):=u(x,t_n)-u_M^n(x),\quad \dot{e}^n(x):=\partial_tu(x,t_n)-\dot{u}_M^n(x),\quad x\in\overline{\Omega},\ 0\leq n\leq\frac{T}{\tau},
\end{equation}
then we have
\begin{theorem}[Error bound of 4th-GIFS]\label{main thm}
Under the assumption (\ref{assumption}), there exist two constants $0<h_0\leq1$ and $0<\tau_0\leq1$
sufficiently small and independent of $\eps$, such that when $0<\tau\leq\tau_0\cdot\min\{\eps^2,h\eps\}$
 and $0<h\leq h_0$, we have
\begin{align}
&\left\|e^n\right\|_{H^1}+\eps^2\left\|\dot{e}^n\right\|_{H^1}\lesssim h^{m_0}+\frac{\tau^4}{\eps^8}, \label{error bound}\\
&\left\|u^n_M\right\|_{L^\infty}\leq C_0+1,\quad \left\|\dot{u}^n_M\right\|_{L^\infty}\leq \frac{C_0+1}{\eps^2},\quad 0\leq n\leq\frac{T}{\tau}.
\label{sol bound}
\end{align}
\end{theorem}
With $u_j^n,\dot{u}_j^n$ obtained from the 4th-GIFP (\ref{4th GIFP s})-(\ref{4th GIFP e}), define the error functions as
\begin{equation*}
e^n(x):=u(x,t_n)-I_M(u^n)(x),\quad \dot{e}^n(x):=\partial_tu(x,t_n)-I_M(\dot{u}^n)(x),\quad x\in\bar{\Omega},\ 0\leq n\leq\frac{T}{\tau},
\end{equation*}
then similarly we have
\begin{theorem}[Error bound of 4th-GIFP]\label{thm Fp}
Under the assumption (\ref{assumption}), there exist two constants $0<h_0\leq1$ and $0<\tau_0\leq1$
sufficiently small and independent of $\eps$, such that when $0<\tau\leq\tau_0\cdot\min\{\eps^2,h\eps\}$
 and $0<h\leq h_0$, we have
\begin{align}
&\left\|e^n\right\|_{H^1}+\eps^2\left\|\dot{e}^n\right\|_{H^1}\lesssim h^{m_0}+\frac{\tau^4}{\eps^8}, \\
&\left\|u^n\right\|_{l^\infty}\leq C_0+1,\quad \left\|\dot{u}^n\right\|_{l^\infty}\leq \frac{C_0+1}{\eps^2},\quad 0\leq n\leq\frac{T}{\tau}.
\end{align}
\end{theorem}

\begin{remark}\label{rm}
In Theorem \ref{main thm} and \ref{thm Fp}, the requirement $0<\tau\leq\tau_0\min\{\eps^2,h\eps\}$ implies that for $\eps=O(1)$, the CFL condition or stability condition is just $\tau\lesssim h$, while for $0<\eps\ll1$, it is $\tau\lesssim\eps^2$ due to essential wave length.
\end{remark}

\subsection{Proof of main result}
For the 4th order Gautschi-type EWI spectral method, in fact, the 4th-GIFS (\ref{4th GIFS s})-(\ref{4th GIFS e}) is a semi-discretization to the KGE,  while the 4th-GIFP (\ref{4th GIFP s})-(\ref{4th GIFP e}) is a full-discretization. For simplicity, we prove the error estimate of the 4th-GIFS, and omit that of the 4th-GIFP which can be done in the same spirit with additional help of interpolation techniques \cite{Cai,Zhao}.

To proceed to the proof of the main result Theorem \ref{main thm}, we first define the projected error
$$e^n_M(x):=P_Mu(x,t_n)-u_M^n(x),\quad \dot{e}^n_M(x):=P_M(\partial_tu(x,t_n))-\dot{u}_M^n(x),\quad 0\leq n\leq\frac{T}{\tau}.$$
Then by triangle inequality and estimates on projection error in \cite{book-sp-GO,Shen} under assumption (\ref{assumption}), we have
\begin{align}
\left\|e^n\right\|_{H^1}+\eps^2\left\|\dot{e}^n\right\|_{H^1}&\lesssim
\left\|e_M^n\right\|_{H^1}+\eps^2\left\|\dot{e}_M^n\right\|_{H^1}+\left\|u(\cdot,t_n)-P_Mu(\cdot,t_n)\right\|_{H^1}\nonumber\\
&\quad+\eps^2\left\|\partial_tu(\cdot,t_n)-P_M(\partial_tu(\cdot,t_n))\right\|_{H^1}\nonumber\\
&\lesssim\left\|e_M^n\right\|_{H^1}+\eps^2\left\|\dot{e}_M^n\right\|_{H^1}+h^{m_0}.\label{proof eq-3}
\end{align}
Thus to prove (\ref{error bound}) in Theorem \ref{main thm}, it is sufficient to work out the corresponding estimate for $e^n_M$ and $\dot{e}^n_M$. The main proof is by the energy method and carried out in the framework of mathematical induction in order to guarantee the boundedness of the numerical solutions \cite{Dong,ZhaoMTI,BDZ,Cai,ZhaoCicp}.
Then the proof is done by the following steps.

Proof of Theorem \ref{main thm}:
For $n=0,$ from the choice of initial data in the scheme, we have
$$e^0=0,\quad \dot{e}^0=0,$$
and results (\ref{error bound}) and (\ref{sol bound}) are obviously true.

For $n\geq1,$
define local truncation errors $\xi_l^n$ and $\dot{\xi}_l^n\,(1\leq n\leq T/\tau,\ l=-M/2,\\ \ldots,M/2-1)$ according to (\ref{4th GIFS m}) as
\begin{subequations}\label{local def}
\begin{align}
\widehat{\xi^{n+1}_l}:=&\widehat{u_l}(t_{n+1})+\widehat{u_l}(t_{n-1})-2\cos(\omega_l\tau)\widehat{u_l}(t_{n})+
\frac{2-2\cos(\omega_l\tau)}{\eps^2\omega_l^2}\widehat{(f(u))}_l(t_n)\nonumber\\
&+\frac{\omega_l^2\tau^2+2\cos(\omega_l\tau)-2}{\eps^2\omega_l^4}\frac{\ud^2}{\ud s^2}\widehat{(f(u))}_l(t_n),\quad 1\leq n\leq \frac{T}{\tau}-1,\\
\widehat{\dot{\xi}^{n+1}_l}:=&\widehat{u_l}'(t_{n+1})-\widehat{u_l}'(t_{n-1})+2\omega_l\sin(\omega_l\tau)\widehat{u_l}(t_{n})+
\frac{2\sin(\omega_l\tau)}{\eps^2\omega_l}\widehat{(f(u))}_l(t_n)\nonumber\\
&+\frac{2\omega_l\tau -2\sin(\omega_l\tau)}{\eps^2\omega_l^3}\frac{\ud^2}{\ud s^2}\widehat{(f(u))}_l(t_n),\quad 1\leq n\leq \frac{T}{\tau}-1,
\end{align}
\end{subequations}
and
\begin{subequations}\label{local def1}
\begin{align}
\widehat{\xi^1_l}:=&\widehat{u_l}(t_{1})-\cos(\omega_l\tau)\widehat{(\phi_1)}_l-\frac{\sin(\omega_l\tau)}{\eps^2\omega_l}\widehat{(\phi_2)}_l+
\frac{1-\cos(\omega_l\tau)}{\eps^2\omega_l^2}\widehat{(f(u))}_l(0)\\
&+\frac{\omega_l\tau-\sin(\omega_l\tau)}{\eps^2\omega_l^3}\frac{\ud}{\ud s}\widehat{(f(u))}_l(0)+\frac{\omega_l^2\tau^2+2\cos(\omega_l\tau)- 2}{2\eps^2\omega_l^4}\frac{\ud^2}{\ud s^2}\widehat{(f(u))}_l(0),\nonumber\\
\widehat{\dot{\xi}^1_l}:=&\widehat{u_l}'(t_{1})+\omega_l\sin(\omega_l\tau)\widehat{(\phi_1)}_l-\frac{\cos(\omega_l\tau)}{\eps^2}\widehat{(\phi_2)}_l+
\frac{\sin(\omega_l\tau)}{\eps^2\omega_l}\widehat{(f(u))}_l(0)\\
&+\frac{1-\cos(\omega_l\tau)}{\eps^2\omega_l^2}\frac{\ud}{\ud s}\widehat{(f(u))}_l(0)+\frac{\omega_l\tau-\sin(\omega_l\tau)}{\eps^2\omega_l^3}\frac{\ud^2}{\ud s^2}\widehat{(f(u))}_l(0).\nonumber
\end{align}
\end{subequations}
%\begin{lemma}[Estimates on $\xi^n$ and $\dot{\xi}^n$]\label{lm:local_error}
%Under the assumption (A), we have estimates for $0<\eps\le 1$
%\begin{equation}
%\left\|\xi^{n}\right\|_{H^1}+\left\|\dot{\xi}^{n}\right\|_{H^1}\lesssim \frac{\tau^4}{\eps^8},\qquad
%n=1,\ldots,\frac{T}{\tau}.
%  \label{lm:local_error_L result}
%\end{equation}
%\end{lemma}
%\begin{proof}

\emph{Step 1}: Estimates on local errors $\xi_l^n$ and $\dot{\xi}_l^n$.

For the solution of the KGE (\ref{KG trun}), $\displaystyle u(x,t_n+s)=\sum_{l=-\infty}^\infty\widehat{u_l}(t_n+s)\phi_l(x)$. So in the Fourier frequency space, we have
$$\eps^2\widehat{u_l}''(t_n+s)+\left(\mu_l^2+\frac{1}{\eps^2}\right)\widehat{u_l}(t_n+s)+\widehat{(f(u))_l}(t_n+s)=0,\quad n\geq0.$$
Then by using the variation-of-constant formula similarly as (\ref{vcf s})-(\ref{vcf}), and noticing (\ref{coefficient}) and (\ref{coefficient0}), we find
\begin{subequations}
\begin{align*}
\widehat{\xi^{n+1}_l}=&-\int_0^\tau\frac{\sin(\omega_l(\tau-w))}{\eps^2\omega_l}\bigg[\widehat{(f(u))_l}(t_n+w)+\widehat{(f(u))_l}(t_n-w)
-2\widehat{(f(u))_l}(t_n)\nonumber\\
&-w^2\frac{\ud^2}{\ud s^2}\widehat{(f(u))_l}(t_n)\bigg]
\ud w,\\
\widehat{\dot{\xi}^{n+1}_l}=&-\int_0^\tau\frac{\cos(\omega_l(\tau-w))}{\eps^2}\bigg[\widehat{(f(u))_l}(t_n+w)+\widehat{(f(u))_l}(t_n-w)
-2\widehat{(f(u))_l}(t_n)\nonumber\\
&-w^2\frac{\ud^2}{\ud s^2}\widehat{(f(u))_l}(t_n)\bigg]
\ud w,\qquad 1\leq n\leq \frac{T}{\tau}-1;\\
\widehat{\xi^{1}_l}\ =&-\int_0^\tau\frac{\sin(\omega_l(\tau-w))}{\eps^2\omega_l}\bigg[\widehat{(f(u))_l}(w)
-\widehat{(f(u))_l}(0)-w\frac{\ud}{\ud s}\widehat{(f(u))_l}(0)\nonumber\\
&-\frac{w^2}{2}\frac{\ud^2}{\ud s^2}\widehat{(f(u))_l}(0)\bigg]
\ud w,\\
\widehat{\dot{\xi}^{1}_l}\ =&-\int_0^\tau\frac{\cos(\omega_l(\tau-w))}{\eps^2}\bigg[\widehat{(f(u))_l}(w)
-\widehat{(f(u))_l}(0)-w\frac{\ud}{\ud s}\widehat{(f(u))_l}(0)\nonumber\\
&-\frac{w^2}{2}\frac{\ud^2}{\ud s^2}\widehat{(f(u))_l}(0)\bigg]
\ud w.
\end{align*}
\end{subequations}
Applying the Taylor's expansion with integral form of the remainder, we get
\begin{subequations}
\begin{align*}
\widehat{\xi^{n+1}_l}=&-\int_0^\tau\frac{\sin(\omega_l(\tau-w))}{\eps^2\omega_l}\frac{w^4}{6}\bigg[\int_0^1
(1-\rho)\bigg(\frac{\ud^4}{\ud s^4}\widehat{(f(u))_l}(t_n+\rho w)\\
&+\frac{\ud^4}{\ud s^4}\widehat{(f(u))_l}(t_n-\rho w)\bigg)\ud \rho\bigg]
\ud w,\\
\widehat{\dot{\xi}^{n+1}_l}=&-\int_0^\tau\frac{\cos(\omega_l(\tau-w))}{\eps^2}\frac{w^4}{6}\bigg[\int_0^1
(1-\rho)\bigg(\frac{\ud^4}{\ud s^4}\widehat{(f(u))_l}(t_n+\rho w)\\
&+\frac{\ud^4}{\ud s^4}\widehat{(f(u))_l}(t_n-\rho w)\bigg)\ud \rho\bigg]
\ud w,\qquad 1\leq n\leq \frac{T}{\tau}-1;\\
\widehat{\xi^{1}_l}\ =&-\int_0^\tau\frac{\sin(\omega_l(\tau-w))}{\eps^2\omega_l}\frac{w^3}{2}\left[\int_0^1
(1-\rho)\frac{\ud^3}{\ud s^3}\widehat{(f(u))_l}(\rho w)\ud\rho\right]\ud w,\\
\widehat{\dot{\xi}^{1}_l}=&-\int_0^\tau\frac{\cos(\omega_l(\tau-w))}{\eps^2}\frac{w^3}{2}\left[\int_0^1
(1-\rho)\frac{\ud^3}{\ud s^3}\widehat{(f(u))_l}(\rho w)\ud\rho\right]\ud w.
\end{align*}
\end{subequations}
Then we have the estimates on the local errors as
\begin{numcases}
\ \left|\widehat{\xi^{1}_l}\right|\ \lesssim\frac{\tau^3}{\sqrt{1+\eps^2\mu_l^2}}\int_0^\tau|\sin(\omega_l(\tau-w))|\int_0^1
\left|\frac{\ud^3}{\ud s^3}\widehat{(f(u))_l}(\rho w)\right|\ud\rho\ud w,\nonumber\\
\left|\widehat{\xi^{n+1}_l}\right|\lesssim\frac{\tau^4}{\sqrt{1+\eps^2\mu_l^2}}\int_0^\tau\left|\sin(\omega_l(\tau-w))\right|\int_0^1\bigg[
\left|\frac{\ud^4}{\ud s^4}\widehat{(f(u))_l}(t_n+\rho w)\right|\nonumber \\
\qquad\qquad+\left|\frac{\ud^4}{\ud s^4}\widehat{(f(u))_l}(t_n-\rho w)\right|\bigg]\ud \rho
\ud w,\quad 1\leq n\leq \frac{T}{\tau}-1,\label{xi est0}
\end{numcases}
and
\begin{numcases}
\ \left|\widehat{\dot{\xi}^{1}_l}\right|\lesssim\frac{\tau^3}{\eps^2}\int_0^\tau\int_0^1
\left|\frac{\ud^3}{\ud s^3}\widehat{(f(u))_l}(\rho w)\right|\ud\rho\ud w,\quad 1\leq n\leq \frac{T}{\tau}-1,\label{dxi est}\\
\left|\widehat{\dot{\xi}^{n+1}_l}\right|\lesssim\frac{\tau^4}{\eps^2}\int_0^\tau\int_0^1
\bigg[\left|\frac{\ud^4}{\ud s^4}\widehat{(f(u))_l}(t_n+\rho w)\right|\nonumber\\
\qquad\qquad+\left|\frac{\ud^4}{\ud s^4}\widehat{(f(u))_l}(t_n-\rho w)\right|\bigg]\ud \rho
\ud w.\nonumber
\end{numcases}
Under condition $\tau\leq\frac{\pi h\eps}{2\sqrt{h^2+4\pi^2\eps^2}}$, which is provided by $\tau\lesssim\min\{\eps^2,h\eps\}$, we have $|\omega_l|\tau\leq\frac{\pi}{2}$ for all $l=-M/2,\ldots,M/2-1$. Then from (\ref{xi est0}) we further have
\begin{numcases}
\ \left|\frac{\widehat{\xi^{1}_l}}{\sin(\omega_l\tau)}\right|\ \lesssim\frac{\tau^3}{\sqrt{1+\eps^2\mu_l^2}}\int_0^\tau\int_0^1
\left|\frac{\ud^3}{\ud s^3}\widehat{(f(u))_l}(\rho w)\right|\ud\rho\ud w,\nonumber\\
\left|\frac{\widehat{\xi^{n+1}_l}}{\sin(\omega_l\tau)}\right|\lesssim\frac{\tau^4}{\sqrt{1+\eps^2\mu_l^2}}\bigg[\int_0^\tau\int_0^1
\bigg|\frac{\ud^4}{\ud s^4}\widehat{(f(u))_l}(t_n+\rho w)\bigg|\ud\rho\ud w\label{xi est}\\
\qquad\qquad\quad\ +\int_0^\tau\int_0^1\bigg|\frac{\ud^4}{\ud s^4}\widehat{(f(u))_l}(t_n-\rho w)\bigg|\ud \rho\ud w\bigg],\quad 1\leq n\leq \frac{T}{\tau}-1.\nonumber
\end{numcases}
With estimates (\ref{dxi est}) and (\ref{xi est}), defining local truncation error functions as
\begin{equation*}
\xi^n(x)=\sum_{l=-N/2}^{N/2-1}\frac{\widehat{\xi^n_l}}{\sin(\omega_l\tau)}\;e^{i\mu_l (x-a)}, \quad
\dot{\xi}^n(x)=\sum_{l=-N/2}^{N/2-1}\widehat{\dot{\xi}^n_l}\;e^{i\mu_l (x-a)},
\quad 1\leq n\leq\frac{T}{\tau},
\end{equation*}
then combining with Paserval's identity and Schwarz's inequality, we get
\begin{align*}
&\|\partial_x\xi^1\|_{H^1}^2+\frac{1}{\eps^2}\|\xi^1\|_{H^1}^2\lesssim\frac{\tau^7}{\eps^2}\int_0^\tau\int_0^1
\left\|\partial_{s}^3f(u)(\cdot,\rho w)\right\|_{H^1}^2\ud\rho\ud w,\\
&\|\partial_x\xi^{n+1}\|_{H^1}^2+\frac{1}{\eps^2}\|\xi^{n+1}\|_{H^1}^2\lesssim\frac{\tau^9}{\eps^2}\int_0^\tau\int_0^1
\left\|\partial_{s}^4f(u)(\cdot,t_n+\rho w)\right\|_{H^1}^2\ud\rho\ud w,\\
&\left\|\dot{\xi}^1\right\|_{H^1}^2\lesssim\frac{\tau^7}{\eps^4}\int_0^\tau\int_0^1
\left\|\partial_{s}^3f(u)(\cdot,\rho w)\right\|_{H^1}^2\ud\rho\ud w,\\
&\left\|\dot{\xi}^{n+1}\right\|_{H^1}^2\lesssim\frac{\tau^9}{\eps^4}\int_0^\tau\int_0^1
\left\|\partial_{s}^4f(u)(\cdot,t_n+\rho w)\right\|_{H^1}^2\ud\rho\ud w,\quad 1\leq n\leq\frac{T}{\tau}-1.
\end{align*}
Thus under assumption (\ref{assumption}), we have for
\begin{subequations}\label{xi bound}
\begin{align}
&\eps^2\left\|\dot{\xi}^1\right\|_{H^1}^2+\|\partial_x\xi^1\|_{H^1}^2+\frac{1}{\eps^2}\|\xi^1\|_{H^1}^2\lesssim\frac{\tau^8}{\eps^{14}},\\
&\eps^2\left\|\dot{\xi}^{n+1}\right\|_{H^1}^2+\|\partial_x\xi^{n+1}\|_{H^1}^2+\frac{1}{\eps^2}\|\xi^{n+1}\|_{H^1}^2
\lesssim\frac{\tau^{10}}{\eps^{18}},\quad1\leq n\leq\frac{T}{\tau}-1.
\end{align}
\end{subequations}

Subtracting the scheme (\ref{4th GIFS e}) from (\ref{local def1}), we get
\begin{align*}
\widehat{e^{1}_l}=\widehat{\xi^1_l},\quad\widehat{\dot{e}^{1}_l}=\widehat{\dot{\xi}^1_l},
\end{align*}
which together with (\ref{xi bound}) and (\ref{proof eq-3}) imply
\begin{align*}
\|e^{1}\|_{H^1}+\eps^2\|\dot{e}^{1}\|_{H^1}\lesssim\frac{\tau^4}{\eps^6}+h^{m_0}.
\end{align*}
Then by triangle inequality and Sobolev's inequality, when $\tau\leq \tau_1\cdot\eps^2$ and $h\leq h_1$,
$$\|u_M^{1}\|_{L^\infty}\leq \|e^{1}\|_{L^\infty}+C_0\leq 1+C_0,\quad \|\dot{u}_M^{1}\|_{L^\infty}\leq \|\dot{e}^{1}\|_{L^\infty}+\frac{C_0}{\eps^2}\leq \frac{1+C_0}{\eps^2},$$
for some constants $\tau_1>0$ and $h_1>0$ independent of $\eps$. Thus (\ref{error bound}) and (\ref{sol bound}) are true for $n=1$.

Now for $n\geq2$, assume (\ref{error bound}) and (\ref{sol bound}) are true for all $n\leq m\leq\frac{T}{\tau}-1$, and then we need to show results (\ref{error bound}) and (\ref{sol bound}) are still valid for $n=m+1$.
Subtracting the scheme (\ref{4th GIFS m}) from (\ref{local def}), we get
\begin{subequations}\label{error eq}
\begin{align}
&\widehat{e^{n+1}_l}+\widehat{e^{n-1}_l}=2\cos(\omega_l\tau)\widehat{e^n_l}+\chi_l^{n+1},\label{error eq a}\\
&\widehat{\dot{e}^{n+1}_l}-\widehat{\dot{e}^{n-1}_l}=-2\omega_l\sin(\omega_l\tau)\widehat{e^n_l}+\dot{\chi}_l^{n+1},\quad 1\leq n\leq \frac{T}{\tau}-1,\label{error eq b}
\end{align}
\end{subequations}
where
\begin{align}
\chi_l^{n+1}:=\widehat{\xi_l^{n+1}}+\widehat{\eta_l^{n+1}},\qquad\dot{\chi}_l^{n+1}:=\widehat{\dot{\xi}_l^{n+1}}+\widehat{\dot{\eta}_l^{n+1}},
\end{align}
with the errors from nonlinear terms as
\begin{subequations}
\begin{align*}
&\widehat{\eta}^{n+1}_l:=\frac{2-2\cos(\omega_l\tau)}{\eps^2\omega_l^2}
\left(\widehat{f^n_l}-\widehat{f_l}(t_n)\right)+\frac{\omega_l^2\tau^2+2\cos(\omega_l\tau)-2}{\eps^2\omega_l^4}
\left(\widehat{\ddot{f}^n_l}-\frac{\ud^2}{\ud s^2}\widehat{f_l}(t_n)\right),\\
&\widehat{\dot{\eta}}^{n+1}_l:=\frac{2\sin(\omega_l\tau)}{\eps^2\omega_l}
\left(\widehat{f^n_l}-\widehat{f_l}(t_n)\right)+\frac{2\omega_l\tau-2\sin(\omega_l\tau)}{\eps^2\omega_l^3}
\left(\widehat{\ddot{f}^n_l}-\frac{\ud^2}{\ud s^2}\widehat{f_l}(t_n)\right).
\end{align*}
\end{subequations}

\emph{Step 2}: Estimates on nonlinear errors $\eta^n_l$ and $\dot{\eta}^n_l$.

Again noting (\ref{coefficient}) with $|\omega_l|\tau\leq\frac{\pi}{2}$, we find
\begin{numcases}
\ \left|\frac{\widehat{\eta}^{n+1}_l}{\sin(\omega_l\tau)}\right|\lesssim\frac{\tau}{\sqrt{1+\eps^2\mu_l^2}}\left|\widehat{f^n_l}-\widehat{f_l}(t_n)\right|
+\frac{\tau^3}{\sqrt{1+\eps^2\mu_l^2}}\left|\widehat{\ddot{f}^n_l}-\frac{\ud^2}{\ud s^2}\widehat{f_l}(t_n)\right|,\nonumber\\
\left|\widehat{\dot{\eta}}^{n+1}_l\right|\lesssim\frac{\tau}{\eps^2}
\left|\widehat{f^n_l}-\widehat{f_l}(t_n)\right|+\frac{\tau^3}{\eps^2}
\left|\widehat{\ddot{f}^n_l}-\frac{\ud^2}{\ud s^2}\widehat{f_l}(t_n)\right|,
\quad 1\leq n\leq \frac{T}{\tau}-1.\label{eta def}
\end{numcases}
Defining nonlinear error functions as
\begin{equation*}
\eta^{n+1}(x)=\sum_{l=-N/2}^{N/2-1}\frac{\widehat{\eta^{n+1}_l}}{\sin(\omega_l\tau)}\;e^{i\mu_l (x-a)}, \
\dot{\eta}^{n+1}(x)=\sum_{l=-N/2}^{N/2-1}\widehat{\dot{\eta}^{n+1}_l}\;e^{i\mu_l (x-a)},
\  1\leq n\leq\frac{T}{\tau}-1,
\end{equation*}
by Parseval's identity and H\"{o}lder's inequality, we have
\begin{numcases}
\ \|\partial_x\eta^{n+1}\|_{H^1}^2+\frac{1}{\eps^2}\|\eta^{n+1}\|_{H^1}^2\lesssim\frac{\tau^2}{\eps^2}\left\|f^n-f(\cdot,t_n)\right\|_{H^1}^2\nonumber\\
\qquad\qquad\qquad\qquad\qquad\qquad\quad+\frac{\tau^6}{\eps^2}\left\|\ddot{f}^n-\partial_s^2f(\cdot,t_n)\right\|_{H^1}^2, \nonumber\\
\eps^2\|\dot{\eta}^{n+1}\|_{H^1}^2\lesssim\frac{\tau^2}{\eps^2}\left\|f^n-f(\cdot,t_n)\right\|_{H^1}^2
+\frac{\tau^6}{\eps^2}\left\|\ddot{f}^n-\partial_s^2f(\cdot,t_n)\right\|_{H^1}^2.\label{proof eq-2}
\end{numcases}
Then under the inductions, for $1\leq n\leq m$, we have
\begin{align}\label{proof eq-1}
\left\|f^n-f(\cdot,t_n)\right\|_{H^1}&=\left\|\int_0^1f'\left(\rho u^n_M+(1-\rho)u(\cdot,t_n)\right)
\ud\rho\left(u^n_M-u(\cdot,t_n)\right)\right\|_{H^1}\nonumber\\
&\lesssim\left\|u^n_M-u(\cdot,t_n)\right\|_{H^1}\lesssim\|e_M^n\|_{H^1}+h^{m_0},
\end{align}
and
\begin{align}
\left\|\ddot{f}^n-\partial_s^2f(\cdot,t_n)\right\|_{H^1}\lesssim&
\left\|f''(u_M^n)\cdot\left(\dot{u}_M^n\right)^2-f''(u(\cdot,t_n))\cdot\left(\partial_su(\cdot,t_n)\right)^2\right\|_{H^1}\nonumber\\
&+\frac{1}{\eps^2}\left\|f'(u_M^n)\partial_{xx}u_M^n-f'(u(\cdot,t_n))\partial_{xx}u(\cdot,t_n)\right\|_{H^1}\nonumber\\
&+\frac{1}{\eps^4}\left\|f'(u_M^n)u_M^n-f'(u(\cdot,t_n))u(\cdot,t_n)\right\|_{H^1}\nonumber\\
&+\frac{1}{\eps^2}\left\|f'(u_M^n)f(u_M^n)-f'(u(\cdot,t_n))f(u(\cdot,t_n))\right\|_{H^1}.\label{proof eq0}
\end{align}
By triangle inequality and Sobolev's inequality, we find
\begin{align}\label{proof eq1}
&\quad\left\|f'(u_M^n)u_M^n-f'(u(\cdot,t_n))u(\cdot,t_n)\right\|_{H^1}\nonumber\\
&\leq
\left\|f'(u_M^n)(u_M^n-u(\cdot,t_n))\right\|_{H^1}+\left\|\left(f'(u_M^n)-f'(u(\cdot,t_n))\right)u(\cdot,t_n)\right\|_{H^1}\nonumber\\
&\lesssim\left\|u_M^n-u(\cdot,t_n)\right\|_{H^1}\lesssim\left\|e_M^n\right\|_{H^1}+h^{m_0},\quad 1\leq n\leq m.
\end{align}
Similarly,
\begin{subequations}\label{proof eq2}
\begin{align}
&\quad\left\|f''(u_M^n)\cdot\left(\dot{u}_M^n\right)^2-f''(u(\cdot,t_n))\cdot\left(\partial_su(\cdot,t_n)\right)^2\right\|_{H^1}\nonumber\\
&\lesssim\left\|\left(\dot{u}_M^n\right)^2-\left(\partial_su(\cdot,t_n)\right)^2\right\|_{H^1}+\frac{1}{\eps^4}
\left\|u(\cdot,t_n)-u_M^n\right\|_{H^1}\nonumber\\
&\lesssim\frac{1}{\eps^2}\left\|\dot{e}_M^n\right\|_{H^1}+\frac{1}{\eps^4}\left\|e_M^n\right\|_{H^1}+\frac{h^{m_0}}{\eps^4},\\
&\quad\left\|f'(u_M^n)f(u_M^n)-f'(u(\cdot,t_n))f(u(\cdot,t_n))\right\|_{H^1}\nonumber\\
&\lesssim\left\|f(u_M^n)-f(u(\cdot,t_n))\right\|_{H^1}+\left\|u_M^n-u(\cdot,t_n)\right\|_{H^1}\nonumber\\
&\lesssim\|e_M^n\|_{H^1}+h^{m_0},\\
&\quad\left\|f'(u_M^n)\partial_{xx}u_M^n-f'(u(\cdot,t_n))\partial_{xx}u(\cdot,t_n)\right\|_{H^1}\nonumber\\
&\lesssim\left\|\partial_{xx}u_M^n-\partial_{xx}u(\cdot,t_n)\right\|_{H^1}+\left\|u_M^n-u(\cdot,t_n)\right\|_{H^1}\nonumber\\
&\lesssim\|e_M^n\|_{H^3}+h^{m_0-2},\qquad 1\leq n\leq m.
\end{align}
\end{subequations}
Plugging (\ref{proof eq1})\&(\ref{proof eq2}) back to (\ref{proof eq0}), under conditions $\tau\lesssim h\eps$ and $\tau\lesssim\eps^2$ we get
\begin{align}
\left\|\ddot{f}^n-\partial_s^2f(\cdot,t_n)\right\|_{H^1}\lesssim&\frac{1}{\eps^2}\left(\left\|\dot{e}_M^n\right\|_{H^1}+\frac{h^{m_0}}{\eps^2}\right)
+\frac{1}{\eps^2}\left(\|e_M^n\|_{H^3}+h^{m_0-2}\right)\nonumber\\
&+\frac{1}{\eps^4}\left(\|e_M^n\|_{H^1}+h^{m_0}\right),\quad 1\leq n\leq m.\label{proof eq3}
\end{align}
Plugging (\ref{proof eq3}) together with (\ref{proof eq-1}) into (\ref{proof eq-2}), under conditions $\tau\lesssim h\eps$ and $\tau\lesssim\eps^2$ we have
\begin{subequations}\label{proof eq4}
\begin{align}
&\|\partial_x\eta^{n+1}\|_{H^1}^2+\frac{1}{\eps^2}\|\eta^{n+1}\|_{H^1}^2\lesssim\tau^2\left(\frac{1}{\eps^2}\|e_M^n\|_{H^1}^2+
\eps^2\left\|\dot{e}_M^n\right\|_{H^1}^2+\frac{h^{2m_0}}{\eps^2}\right), \\
&\eps^2\|\dot{\eta}^{n+1}\|_{H^1}^2\lesssim\tau^2\left(\frac{1}{\eps^2}\|e_M^n\|_{H^1}+\eps^2\left\|\dot{e}_M^n\right\|_{H^1}+\frac{h^{2m_0}}{\eps^2}\right),
\quad 1\leq n\leq m.
\end{align}
\end{subequations}

\emph{Step 3}: Error equations on $e^n_l$ and $\dot{e}^n_l$.

Denote funcation
$$\mathcal{T}_l^n:=\sin(\omega_l(n+1)\tau),\qquad n\geq-1.$$
Multiplying both sides of (\ref{error eq a}) by $\mathcal{T}_l^{m-n}$ and then summing up for $1\leq n\leq m$, we get
\begin{align*}
\quad\sum_{n=1}^m\mathcal{T}_l^{m-n}\left(\widehat{e_l^{n+1}}+\widehat{e_l^{n-1}}\right)&=\sum_{n=1}^m\mathcal{T}_l^{m-n}(2\cos(\omega_l\tau)\widehat{e_l^{n}}+u_l^{n+1})\\
&=\sum_{n=1}^m[(\mathcal{T}_l^{m+1-n}+\mathcal{T}_l^{m-1-n})\widehat{e_l^{n}}+\mathcal{T}_l^{m-n}u_l^{n+1}]\\
&=\sum_{n=0}^{m-1}\mathcal{T}_l^{m-n}\widehat{e_l^{n+1}}+\sum_{n=2}^{m+1}\mathcal{T}_l^{m-n}\widehat{e_l^{n-1}}+\sum_{n=1}^{m}\mathcal{T}_l^{m-n}u_l^{n+1},
\end{align*}
which consequently shows
\begin{equation}\label{proof eq5}
\mathcal{T}_l^{0}\widehat{e_l^{m+1}}=-\mathcal{T}_l^{m-1}\widehat{e_l^{0}}+\mathcal{T}_l^{m}\widehat{e_l^{1}}+\sum_{n=1}^{m}\mathcal{T}_l^{m-n}u_l^{n+1}.
\end{equation}
Noting $e_l^{0}=0$ and by Cauchy's inequality,
\begin{align*}
|\widehat{e_l^{m+1}}|^2&\leq2\left[|\mathcal{T}_l^{m}|^2\left|\frac{\widehat{e_l^{1}}}{\sin(\omega_l\tau)}\right|^2+
m\sum_{n=1}^{m}|\mathcal{T}_l^{m-n}|^2\left|\frac{u_l^{n+1}}{\sin(\omega_l\tau)}\right|^2 \right]\\
&\leq2\left[\left|\frac{\widehat{\xi_l^{1}}}{\sin(\omega_l\tau)}\right|^2+
2m\sum_{n=1}^{m}\left(\left|\frac{\widehat{\xi_l^{n+1}}}{\sin(\omega_l\tau)}\right|^2+
\left|\frac{\widehat{\eta_l^{n+1}}}{\sin(\omega_l\tau)}\right|^2\right)\right].
\end{align*}
Multiplying both sides by $(1+\mu_l^2)\frac{1}{\eps^2}$, %and $(1+\mu_l^2)\mu_l^2$, respectively,
and then summing up for $l=-M/2,\ldots,M/2-1$, by (\ref{xi bound}) and (\ref{proof eq4}), we get
\begin{align}
\frac{1}{\eps^2}\left\|e_M^{m+1}\right\|_{H^1}^2&\lesssim \frac{1}{\eps^2}\|\xi^1\|_{H^1}^2+\frac{m}{\eps^2}
\sum_{n=1}^{m}\left(\left\|\xi^{n+1}\right\|_{H^1}^2+\left\|\eta^{n+1}\right\|_{H^1}^2\right)\nonumber\\
&\lesssim\frac{1}{\eps^2}\left(\frac{\tau^8}{\eps^{16}}+h^{2m_0}\right)+\tau\sum_{n=1}^{m}\left(\frac{1}{\eps^2}\left\|e_M^n\right\|_{H^1}^2
+\eps^2\left\|\dot{e}_M^{n}\right\|_{H^1}^2\right).\label{em}
%\left\|\partial_xe_M^{m+1}\right\|_{H^1}^2&\lesssim \|\partial_x\xi^1\|_{H^1}^2+m
%\sum_{n=1}^{m}\left(\left\|\partial_x\xi^{n+1}\right\|_{H^1}^2+\left\|\partial_x\eta^{n+1}\right\|_{H^1}^2\right)\nonumber\\
%&\lesssim\frac{1}{\eps^2}\left(\frac{\tau^8}{\eps^{16}}+h^{2m_0}\right)+\tau\sum_{n=1}^{m}\left(\frac{1}{\eps^2}\left\|e^n\right\|_{H^1}^2
%+\eps^2\left\|\dot{e}_l^{n}\right\|_{H^1}^2\right).\label{em b}
\end{align}
Plugging (\ref{proof eq5}) into (\ref{error eq b}), we get
$$\widehat{\dot{e}^{m+1}_l}-\widehat{\dot{e}^{m-1}_l}=-2\omega_l\left(\mathcal{T}_l^{m-1}\widehat{e_l^{1}}+\sum_{n=1}^{m-1}\mathcal{T}_l^{m-1-n}\chi_l^{n+1}\right)
+\dot{\chi}_l^{m+1}.$$
Using the recurrence equation, for some odd $m=2p-1(p\geq1)$, we have
\begin{align*}
&\quad\widehat{\dot{e}^{2p}_l}-\widehat{\dot{e}^{0}_l}\\
&=\sum_{n=1}^{p}\dot{\chi}_l^{2n}-2\omega_l\sum_{n=1}^{p}\sin(\omega_l(2n-1)\tau)\widehat{e_l^{1}}
-2\omega_l\sum_{q=2}^{p}\sum_{n=1}^{2q-2}\sin(\omega_l(2q-1-n)\tau)\chi_l^{n+1}\\
&=\sum_{n=1}^{p}\dot{\chi}_l^{2n}-2\omega_l\sin^2(\omega_lp\tau)\frac{\widehat{\xi_l^{1}}}{\sin(\omega_l\tau)}-
2\omega_l\sum_{n=1}^{2p-2}\chi^{n+1}_l\sum_{q=q_n}^{p}\sin(\omega_l(2q-1-n)\tau),
\end{align*}
with $q_n:=\lceil\frac{n}{2}\rceil+1$, where $\lceil\cdot\rceil$ denotes the ceiling function.
Then we have
\begin{align*}
\eps^2\left\|\dot{e}_M^{2p}\right\|_{H^1}^2
\lesssim& \eps^2p\sum_{n=1}^{p}\left\|\dot{\chi}^{2n}\right\|_{H^1}^2+\left\|\partial_x\xi^{1}\right\|_{H^1}^2+\frac{1}{\eps^2}\left\|\xi^{1}\right\|_{H^1}^2+
p\sum_{n=1}^{2p-2}\left\|\partial_x\chi^{n+1}\right\|_{H^1}^2\\
&+p\sum_{n=1}^{2p-2}\frac{1}{\eps^2}\left\|\chi^{n+1}\right\|_{H^1}^2.
\end{align*}
Similarly, we can get an estimate for the case $m=2p(p\geq1)$ as
\begin{align*}
\eps^2\left\|\dot{e}_M^{2p+1}\right\|_{H^1}^2
\lesssim& \eps^2\left\|\dot{e}_M^{1}\right\|_{H^1}^2+\eps^2p\sum_{n=1}^{p}\left\|\dot{\chi}^{2n+1}\right\|_{H^1}^2+\left\|\partial_x\xi^{1}\right\|_{H^1}^2
+\frac{1}{\eps^2}\left\|\xi^{1}\right\|_{H^1}^2\\
&+p\sum_{n=1}^{2p-1}\left\|\partial_x\chi^{n+1}\right\|_{H^1}^2
+p\sum_{n=1}^{2p-1}\frac{1}{\eps^2}\left\|\chi^{n+1}\right\|_{H^1}^2.
\end{align*}
All together, we have
\begin{align*}
\eps^2\left\|\dot{e}_M^{m+1}\right\|_{H^1}^2
\lesssim& \eps^2\left\|\dot{e}_M^{1}\right\|_{H^1}^2+m\sum_{n=1}^{m}\eps^2\left\|\dot{\chi}^{n+1}\right\|_{H^1}^2+\left\|\partial_x\xi^{1}\right\|_{H^1}^2
+\frac{1}{\eps^2}\left\|\xi^{1}\right\|_{H^1}^2\\
&+m\sum_{n=1}^{m}\left(\left\|\partial_x\chi^{n}\right\|_{H^1}^2
+\frac{1}{\eps^2}\left\|\chi^{n}\right\|_{H^1}^2\right).
\end{align*}
Then by (\ref{xi bound}) and (\ref{proof eq4}), we have
\begin{align}\label{dem}
\eps^2\left\|\dot{e}_M^{m+1}\right\|_{H^1}^2
\lesssim \frac{1}{\eps^2}\left(\frac{\tau^8}{\eps^{16}}+h^{2m_0}\right)
+\tau\sum_{n=1}^{m}\left(\eps^2\left\|\dot{e}_M^{n}\right\|_{H^1}^2
+\frac{1}{\eps^2}\left\|e_M^{n}\right\|_{H^1}^2\right).
\end{align}
Adding up (\ref{em}) and (\ref{dem}), we get
\begin{align*}
&\quad\eps^2\left\|\dot{e}_M^{m+1}\right\|_{H^1}^2+
\frac{1}{\eps^2}\left\|e_M^{m+1}\right\|_{H^1}^2\\
&\lesssim
\frac{1}{\eps^2}\left(\frac{\tau^8}{\eps^{16}}+h^{2m_0}\right)+\tau\sum_{n=1}^{m}\left(\eps^2\left\|\dot{e}_M^{n}\right\|_{H^1}^2
+\frac{1}{\eps^2}\left\|e_M^{n}\right\|_{H^1}^2\right),
\end{align*}
and then by discrete Gronwall's inequality, we get
$$\eps^2\left\|\dot{e}_M^{m+1}\right\|_{H^1}+\left\|e_M^{m+1}\right\|_{H^1}\\
\lesssim \frac{\tau^4}{\eps^{8}}+h^{m_0}.$$
Then by triangle inequality and Sobolev's inequality together with (\ref{proof eq-3}), when $\tau\leq \tau_2\cdot\eps^2$ and $h\leq h_2$,
$$\|u_M^{m+1}\|_{L^\infty}\leq \|e^{m+1}\|_{L^\infty}+C_0\leq 1+C_0,\quad \|\dot{u}_M^{m+1}\|_{L^\infty}\leq \|\dot{e}^{m+1}\|_{L^\infty}+\frac{C_0}{\eps^2}\leq \frac{1+C_0}{\eps^2},$$
for some constants $\tau_2>0$ and $h_2>0$ independent of $\eps$. Thus (\ref{error bound}) and (\ref{sol bound}) are true for $n=m+1$, and the proof is completed by choosing $\tau_0=\min\{\tau_1,\tau_2\},$ $ h_0=\min\{h_1,h_2\}$. \qed

\begin{remark}
The proof technique here is different from that in \cite{Dong,BDZ}. The proof used in \cite{Dong,BDZ} can hardly get the rigorous error estimates for the group of trigonometric integrators proposed in \cite{DongNew}, while it is believed that the proof established here could offer some clues to that which will be our future work.
\end{remark}

To close this chapter, we make an important remark on another potential application of the proposed method.
The recent developed multiscale time integrators (MTIs) in \cite{Zhao,ZhaoMTI,ZhaoKGS,CaiDirac} only achieved the first order uniform accuracy for solving the highly-oscillatory equations. All of them are using the second order EWIs as the key integration tools. Now with the higher order EWIs and using higher order multiscale expansion in corresponding context, MTIs with higher order of uniform accuracy could be proposed, which is going to appear in our future work.

\section{Numerical results}\label{sec: result}
In this section, we present the numerical results of the 4th-GIFP (\ref{4th GIFP s})-(\ref{4th GIFP e}) and a 6th order GIFP (shorted as 6th-GIFP) proposed in Section \ref{subsec: EWIs}. As comparisons, we also present the numerical results of the GIFP method proposed in \cite{Dong} (shorted as GIFP) and the classical 4th order Runger-Kutta method \cite{Iserles,Hair} with Fourier spectral discretization (shorted as RK4FP). Throughout the section, we consider the KGE (\ref{KG trun}) with cubic nonlinearity, i.e.
$$f(u)=\lambda u^3,\quad \lambda\in\bR,$$ which occurs in the most application cases and physical situations \cite{Dong,Faou,Zhao,ZhaoMTI,Masmoudi,Machihara1,Machihara2,Duncan}. Choose
$$\lambda=1, \quad \phi_1=2\fe^{-x^2},\quad \phi_2=3\fe^{-x^2},\quad x\in\Omega,$$
in (\ref{KG trun}),
where the `exact' solution of the problem is obtained via the 6th-GIFP method with very small time step and mesh size, e.g. $\tau=1E-5,\, h=1/16$.
We choose the bounded interval $\Omega=[-32,32]$, i.e. $b=-a=32$,
which is large enough to guarantee that the periodic boundary condition does not
 introduce a significant aliasing error relative to the original problem.
To measure the error, we compute the $H^1$-norm of the error
$$e(x,T)=u(x,T)-I_Mu^N(x),\quad \mbox{where}\quad N=\frac{T}{\tau},$$
for some fixed time $T>0$.

Firstly, we shall test the temporal convergence rate of the proposed 4th-GIFP and 6th-GIFP for a fixed $0<\eps<1$ in the normal (relativistic) regime, i.e. $\eps=O(1)$. The numerical results at $T=2$ under different $\tau$ are given in Table \ref{tab: convergence}. As $\eps$ becomes small, we understand the necessity of condition $\tau\lesssim\eps^2$ from either stability or accuracy point of view. However for fixed $\eps=O(1)$ in the normal regime, the stability constraint $\tau\lesssim h$ imposed in Theorem \ref{main thm} and Remark \ref{rm} is mainly used for the rigorous mathematical proof. Thus, in addition we test the error of the methods for solving the KGE with $\eps=0.5$ and a fixed large $\tau=0.1$ but under different small mesh size $h$. The results are shown in Table \ref{tab: stab}.

Then we study the errors and meshing strategy of the 4th-GIFP and 6th-GIFP in the nonrelativistic limit regime, i.e. $0<\eps\ll1$. We test and study the temporal and spatial error separately.
Table. \ref{tb: spa} shows the spatial error of numerical methods at $T=2$
under different $\eps$ and $h$ with a very small time step $\tau=10^{-5}$
such that the discretization error in time is negligible.
Table. \ref{tab: scalabity} shows the temporal error of numerical methods at $T=2$
under different $\eps$ and $\tau$ with a small mesh size $h=1/16$ such that
the discretization error in space is negligible.

\begin{figure}[t!]
\center
\includegraphics[height=6cm,width=13cm]{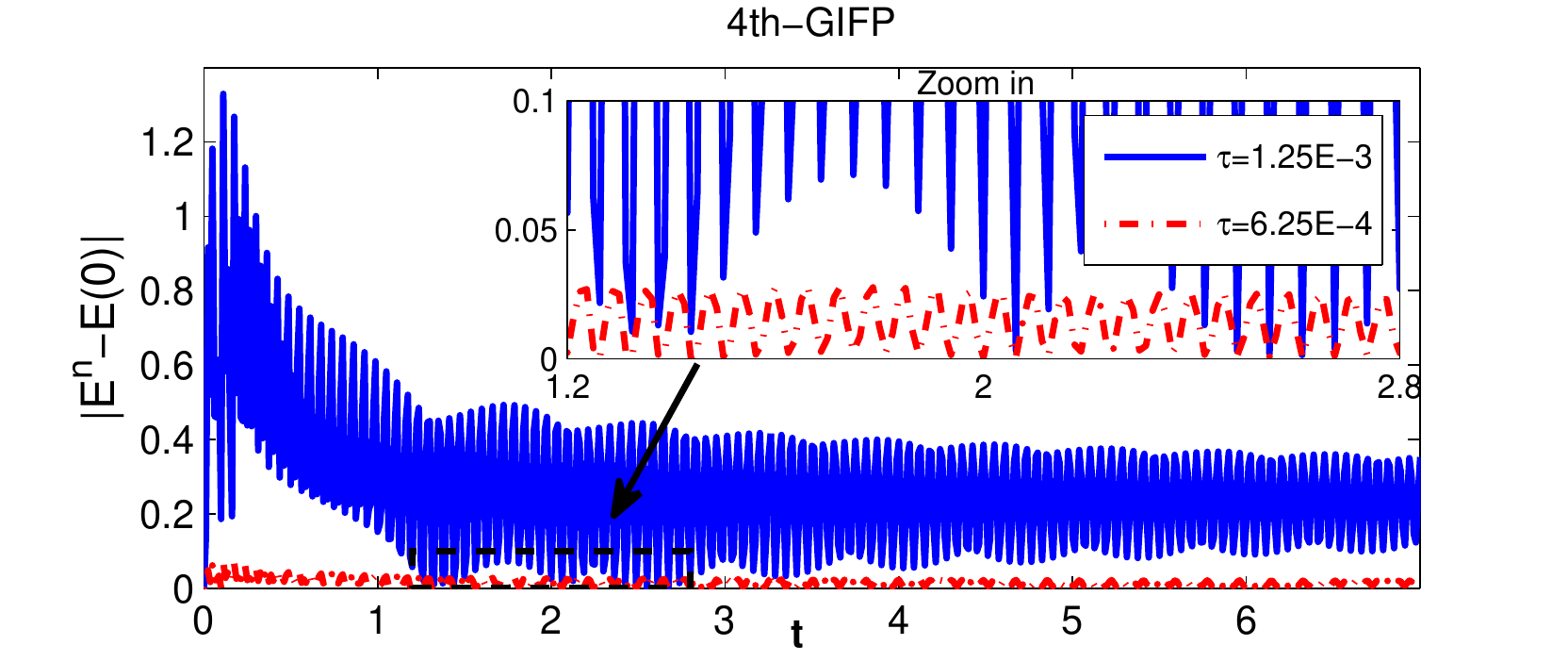}
\includegraphics[height=6cm,width=13cm]{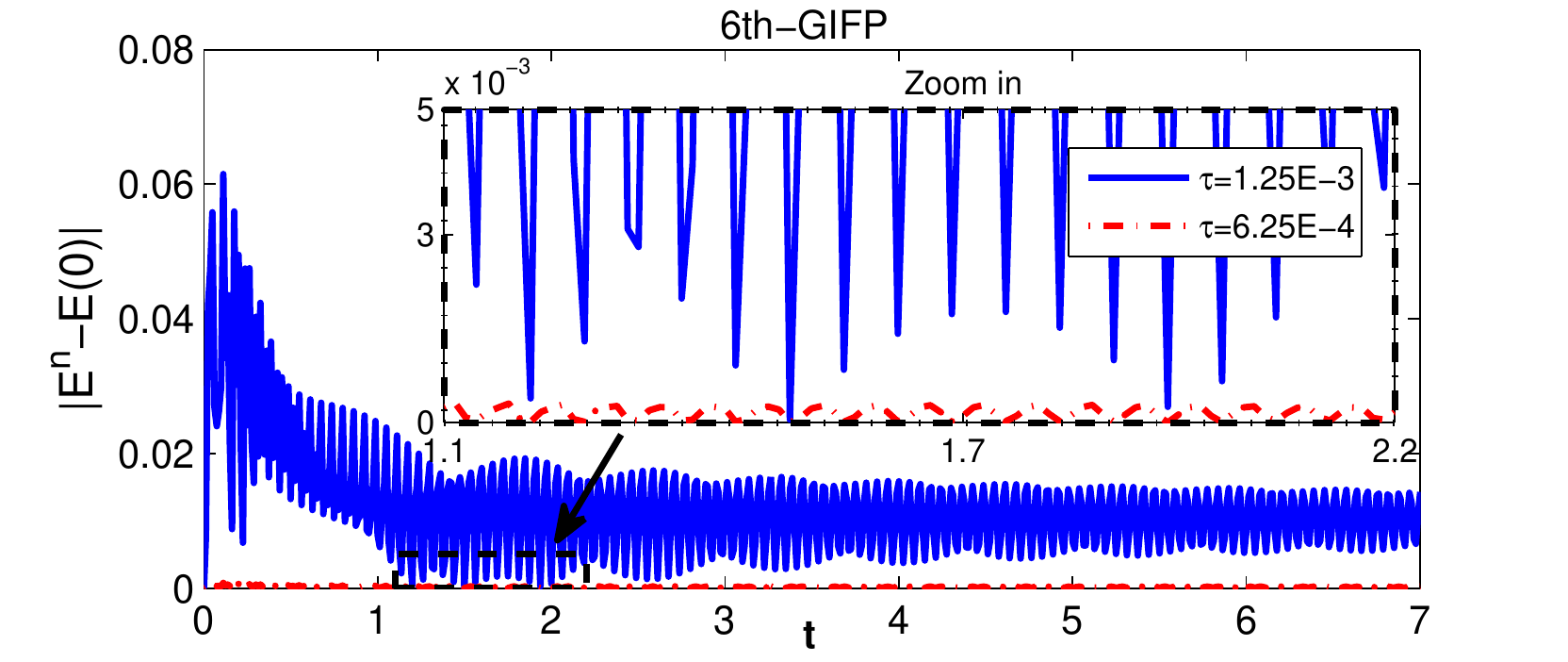}
\includegraphics[height=6cm,width=13cm]{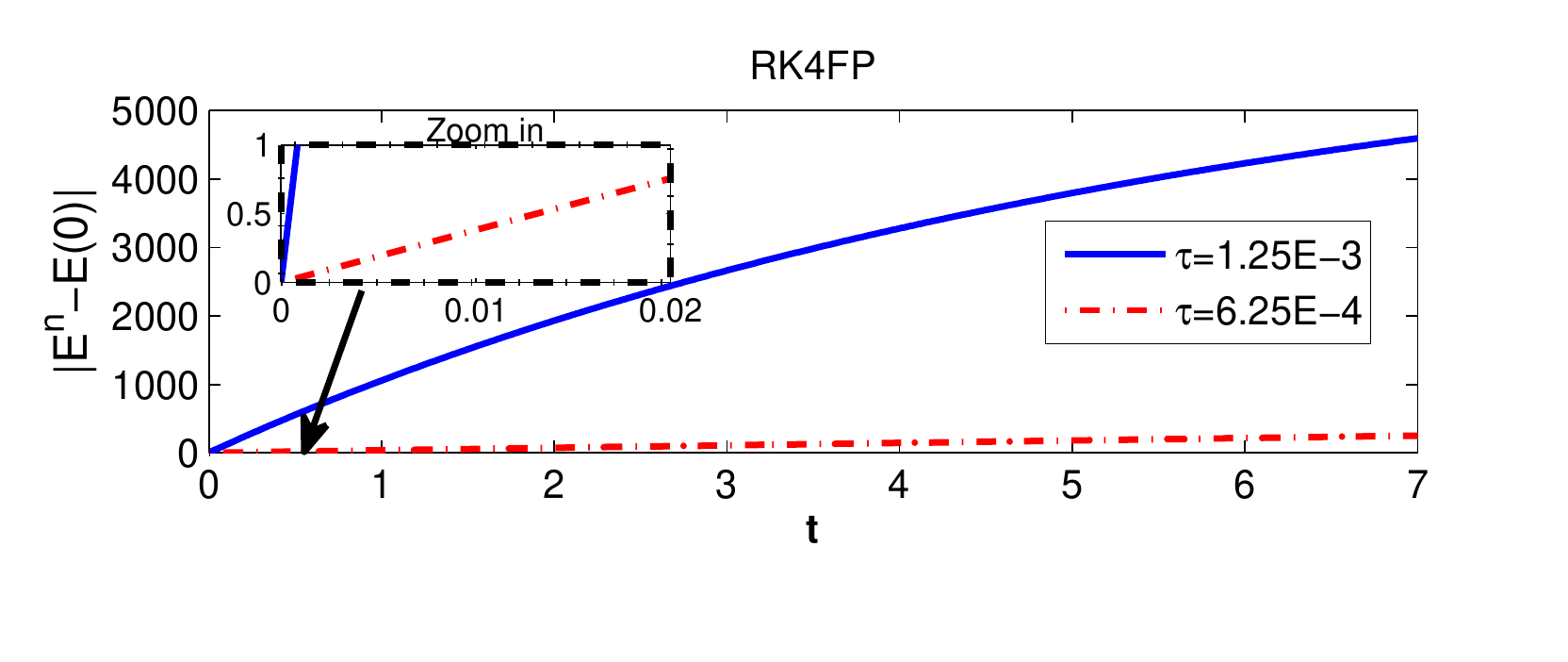}
\caption{Energy error of the 4th-GIFP\&6th-GIFP and comparisons with RK4FP for solving the KGE with $\epsilon=0.05$ under $h=1/8$ and different $\tau$.}\label{fig:0}
\end{figure}

At last, we test the energy conservation property of the schemes. The energy errors of 4th-GIFP and 6th-GIFP during the computation, i.e. the error between the exact energy $E(t)=E(0)$ and the numerical energy
$$E^n:=\int_{a}^b\left[\eps^2|I_M\dot{u}^n(x)|^2+|\partial_xI_Mu^n(x)|^2+\frac{1}{\eps^2}|I_Mu^n(x)|^2+F(I_Mu^n(x))\right]d x,$$
under a small mesh size $h=1/8$, are plotted in Figure \ref{fig:0} together with comparisons with results of the RK4FP.

\begin{table}[t!]
\tabcolsep 0pt
\caption{Temporal error and convergence rate test of the 4th-GIFP\&6th-GIFP with a fixed $\eps=0.5$ in normal regime:
$\|e(\cdot,T)\|_{H^1}$ at $T=2$ for different $\tau$ with $h=1/16$.}\label{tab: convergence}
\begin{center}
\begin{tabular*}{1\textwidth}{@{\extracolsep{\fill}}lllll}
\hline
4th-GIFP                 & $\tau_0=0.1$	   &$\tau_0/2$       &$\tau_0/4$ &	$\tau_0/8$\\ \hline
$\|e(\cdot,T)\|_{H^1}$    &4.55E-02	&1.60E-03	&9.52E-05	&5.90E-06\\
rate	                &--         &4.84	    &4.07	    &4.01\\
\hline
6th-GIFP	                & $\tau_0=0.1$	   &$\tau_0/2$       &$\tau_0/4$ &$\tau_0/8$	 \\ \hline
$\|e(\cdot,T)\|_{H^1}$   &4.50E-03	&2.43E-05	&3.74E-07	&5.74E-09\\
rate	                &--                &7.53	&6.07	&6.02
\\ \hline
\end{tabular*}
\end{center}
\end{table}

\begin{table}[t!]
\tabcolsep 0pt
\caption{Stability test of the 4th-GIFP\&6th-GIFP with a fixed $\varepsilon=0.5$ in normal regime:
$\|e(\cdot,T)\|_{H^1}$ at $T=2$ for different $h$ with $\tau=0.1$.}\label{tab: stab}
\begin{center}
\begin{tabular*}{1\textwidth}{@{\extracolsep{\fill}}llll}
\hline
4th-GIFP                 & $h=1/8$	   &$h=1/16$       &$h=1/32$ 	\\ \hline
$\|e(\cdot,T)\|_{H^1}$     &4.55E-02	&4.55E-02	&4.55E-02	\\
\hline
6th-GIFP	                & $h=1/8$	   &$h=1/16$       &$h=1/32$  \\ \hline
 $\|e(\cdot,T)\|_{H^1}$   &4.50E-03	&4.50E-03	&5.85E-02	
\\ \hline
\end{tabular*}
\end{center}
\end{table}

\begin{table}[t!]
\tabcolsep 0pt
\caption{Spatial error of the 4th-GIFP\&6th-GIFP in nonrelativistic limit regime: $\|e(\cdot,T)\|_{H^1}$ at $T=2$ for different $\eps$ and $h$ with  $\tau=10^{-5}$.}\label{tb: spa}
\begin{center}
\begin{tabular*}{1\textwidth}{@{\extracolsep{\fill}}lllll}
\hline
4th-GIFP         & $h_0=1$	      &$h_0/2$	        &$h_0/4$	        &$h_0/8$	\\
\hline
$\eps_0=0.1$	               &6.89E+00	&7.28E-01	&4.18E-04	&6.18E-08\\
$\eps_0/2$	               &6.94E+00	&1.06E+00	&5.58E-04	&1.71E-08\\
$\eps_0/4$	               &7.34E+00	&1.09E+00	&6.16E-04	&5.15E-09\\
$\eps_0/8$	               &7.21E+00	&1.12E+00	&5.38E-04	&6.62E-07\\
\hline
6th-GIFP         & $h_0=1$	      &$h_0/2$	        &$h_0/4$	        &$h_0/8$	\\
\hline
$\eps_0=0.1$	               &6.89E+00	&7.28E-01	&4.18E-04	&6.19E-08\\
$\eps_0/2$	               &6.94E+00	&1.06E+00	&5.58E-04	&1.72E-08\\
$\eps_0/4$	               &7.34E+00	&1.09E+00	&6.16E-04	&4.54E-09\\
$\eps_0/8$	               &7.21E+00	&1.12E+00	&5.38E-04	&1.32E-09\\
\hline
\end{tabular*}
\end{center}
\end{table}

\begin{table}[t!]
\tabcolsep 0pt
\caption{Temporal error and meshing strategy of the 4th-GIFP\&6th-GIFP with convergence rate and comparisons with GIFP and RK4FP in nonrelativistic limit regime:
$\|e(\cdot,T)\|_{H^1}$ at $T=2$ for different $\eps$ and $\tau$ under $\tau=O(\eps^2)$ with $h=1/16$.}\label{tab: scalabity}
\begin{center}
\begin{tabular*}{1\textwidth}{@{\extracolsep{\fill}}llll}
\hline
4th-GIFP                 & $\tau_0=1.25*10^{-3}$	   &$\tau_0/2$       &$\tau_0/4$	\\ \hline
$\eps_0=0.05,\tau_0$    &3.47E-03	&2.04E-04	&1.26E-05\\
rate	                &--         &4.09	     &4.02\\
$\eps_0/2,\tau_0/2^2$   &2.80E-03	&1.69E-04	&1.05E-05\\
rate	                &--         &4.05	     &4.01\\
$\eps_0/4,\tau_0/4^2$   &2.56E-03	&1.55E-04	&9.65E-06\\
rate	                &--         &4.04	     &4.01\\
\hline
6th-GIFP	                & $\tau_0=1.25*10^{-3}$	   &$\tau_0/2$       &$\tau_0/4$	 \\ \hline
$\eps_0=0.05,\tau_0$    &4.70E-05	&6.27E-07	&9.67E-09\\
rate	                &--                &6.23	      &6.02\\
$\eps_0/2,\tau_0/2^2$   &2.05E-05	&2.89E-07	&4.76E-09\\
rate	                &--                &6.14	      &5.93\\
$\eps_0/4,\tau_0/4^2$   &1.34E-05	&1.98E-07	&3.10E-09\\
rate	                &--                &6.07	      &6.00
\\ \hline
GIFP                  & $\tau_0=1.25*10^{-3}$	   &$\tau_0/2$       &$\tau_0/4$	 \\ \hline
$\eps_0=0.05,\tau_0$    &3.25E-01  &7.94E-02	&1.97E-02\\
rate	                &--                &2.03	      &2.01\\
$\eps_0/2,\tau_0/2^2$   &3.12E-01	&7.60E-02	&1.86E-02\\
rate	                &--                &2.04	      &2.03\\
$\eps_0/4,\tau_0/4^2$   &3.38E-01	&7.31E-02	&1.80E-02\\
rate	                &--                &2.20	      &2.02\\
\hline
RK4FP                  & $\tau_0=1.25*10^{-3}$	   &$\tau_0/2$       &$\tau_0/4$	 \\ \hline
$\eps_0=0.05,\tau_0$    &3.53E+00  &2.21E-01	&1.25E-02\\
rate	                &--                &4.00	      &4.13\\
$\eps_0/2,\tau_0/2^2$   &8.06E+00	&8.32E-01	&4.78E-02\\
rate	                &--                &3.28	      &4.12\\
$\eps_0/4,\tau_0/4^2$   &6.38E+00	&3.03E+00	&1.83E-01\\
rate	                &--                &1.08	      &4.04\\
\hline
\end{tabular*}
\end{center}
\end{table}

Based on results from Tables \ref{tab: convergence}-\ref{tab: scalabity} and Figure \ref{fig:0}, we can draw the following observations:

(i) The 4th-GIFP and 6th-GIFP have 4th and 6th temporal accuracy order, respectively, and both of them have the spectral accuracy in spatial discretization. Our theoretical error bound is optimal. The theoretical stability constrain $\tau\lesssim h$ when $\eps=O(1)$ does not seem to be essential in computing.

(ii) As $\eps$ decreases to zero in the nonrelativistic limit regime, the meshing strategy of the 4th-GIFP and 6th-GIFP is $\tau=O(\eps^2)$ and $h=O(1)$. Under the same meshing strategy, the computational error of them is much smaller than that of the GIFP proposed in \cite{Dong} and the classical RK4FP method.

(iii) The 4th-GIFP and 6th-GIFP conserve the energy very well. The energy obtained from the numerical
solution is just a small fluctuation from the exact energy during the computation, while in contrast, the energy error of the RK4FP keeps growing. As time step $\tau$ decreases to zero, the energy error converges to zero.

\section{Conclusions}\label{sec: con}
A group of high order Gautschi-type exponential wave integrators (EWIs) Fourier pseudospectral method were proposed and analyzed for solving the Klein-Gordon equation (KGE) in the nonrelativistic limit regime with a dimensionless parameter $0<\eps\ll1$, where the small $\eps$ makes the solution of the problem propagates waves with wavelength $O(\eps^2)$ in time axis, i.e. high oscillations occur in time. The scheme is fully explicit and time symmetric. In fact, we proposed a way to construct an EWI spectral method with temporal accuracy at any even order and spectral spatial accuracy, provided the solution of the problem is smooth enough. Rigorous error estimates were established to show the meshing strategy of the proposed methods is $\tau=O(\eps^2)$ and $h=O(1)$, as $0<\eps\ll1$ in the nonrelativistic limit regime, where $\tau$ and $h$ denote the time step and mesh size respectively. In view of the essential wave length propagating in time, the proposed EWIs Fourier pseudospectral method offer the high order convergence rate with the `optimal' meshing strategy among all classical numerical methods for directly solving the KGE in the limit regime. The proposed method also implies a promising way to construct multiscale time integrators \cite{Zhao,ZhaoMTI} with higher order uniform accuracy in future.  Extensive numerical experiments were done to confirm the theoretical accuracy order and the meshing strategy. Comparisons with existing classical methods were carried out to show the superiority of the proposed methods. It is also believed that the proposed methods can find wide applications in effectively solving other KG-type oscillatory equations or coupled systems in future work.

% put your thanks here
\section*{Acknowledgments}
This work is supported by the French ANR project MOONRISE ANR-14-CE23-0007-01.

\end{document}